\documentclass[sii]{ipart}

\usepackage{graphicx}
\RequirePackage[numbers]{natbib}
\usepackage{amsthm}
\usepackage{multirow}
\usepackage[ruled,vlined]{algorithm2e}
\usepackage{float}
\usepackage{breakcites}
\usepackage{epsfig}
\usepackage{amssymb}
\usepackage{amsmath}
\usepackage{hyperref}
\usepackage{verbatim}
\usepackage{subcaption}
\usepackage{tikz}
\usepackage{color}

\usepackage{bbm}
\usetikzlibrary{arrows}
\usepackage[margin=1.2in]{geometry}
\newcommand\given[1][]{\:#1\vert\:}
\DeclareMathOperator{\E}{\mathbb{E}}
\tikzset{cross/.style={cross out, draw=black, minimum size=2*(#1-\pgflinewidth), inner sep=0pt, outer sep=0pt},
cross/.default={5pt}}
\usetikzlibrary{shapes.misc}
\usepackage{array}
\newcommand{\be}{\begin{eqnarray}}
\newcommand{\ee}{\end{eqnarray}}
\newcommand{\bea}{\begin{eqnarray*}}
\newcommand{\eea}{\end{eqnarray*}}

\usepackage{nicefrac}
\numberwithin{equation}{section}
\DeclareMathOperator*{\esssup}{ess\,sup}
 \newcolumntype{P}[1]{>{\centering\arraybackslash}m{#1}}

\newtheorem{lemma}{Lemma}

\startlocaldefs
\endlocaldefs

\firstpage{1}
\lastpage{1}
\begin{document}

%
%
%\author{ Jack Noonan 
%       %etc.
%}
%
%%\authorrunning{Short form of author list} % if too long for running head
%
%\institute{           J. Noonan \at
%              School of Mathematics, Cardiff University, Cardiff, CF24 4AG, UK \\
%               \email{Noonanj1@cardiff.ac.uk}
%% The correct dates will be entered by the editor
%
%
%\maketitle
%
%
%\begin{abstract}
%
%
%
%
%\keywords{}
%% \PACS{PACS code1 \and PACS code2 \and more}
%% \subclass{MSC code1 \and MSC code2 \and more}
%\end{abstract}
%
%

\begin{frontmatter}

%%  "Title of the Paper"
\title{Online change-point detection for a transient change}
%\thankstext{T1}{???}
\begin{aug}
    \author{\inits{J.}\fnms{Jack} \snm{Noonan}\ead[label=e1]{Noonanj1@cardiff.ac.uk}}
    \address{School of Mathematics, Cardiff University, Cardiff, CF24 4AG, UK\\
             \printead{e1}}
%    \author{\inits{S.}\fnms{Second} \snm{Author}\ead[label=e2]{second@somewhere.com}}
%    \address{Address of the Second Author\\
%             Country\\
%             \printead{e2}}
%    \and
%    \author{\inits{T.~N.}\fnms{Third Name} \snm{Author}
%            \ead[label=e3]{third@somewhere.com}%
%            \ead[label=u1,url]{http://www.foo.com}}
%    \address{Address of the Third Author\\
%             Country\\
%             \printead{e3}\\
%             \printead{u1}}
%    \thankstext{t2}{Corresponding author.}
\end{aug}
\received{\sday{6} \smonth{4} \syear{2021}}

\begin{abstract}
We consider a popular online change-point problem of detecting a transient change in distributions of i.i.d. random variables. For this change-point problem, several change-point procedures are formulated and some advanced results for a particular procedure are surveyed. Some new approximations for the average run length to false alarm are offered and the power of these procedures for detecting a transient change in mean of a sequence of normal random variables is compared. 
\end{abstract}

%\begin{keyword}[class=AMS]
%\kwd[Primary ]{}
%\kwd{}
%\kwd[; secondary ]{}
%\end{keyword}

%%  Upper case for every keyword
%\begin{keyword}
%\kwd{}
%\kwd{}
%\end{keyword}

%\tableofcontents


%%  Upper case for every keyword

%\tableofcontents
\end{frontmatter}

\section{Introduction}

The subject of change-point detection (or statistical quality control) is devoted to monitoring and detecting changes in the structure of a time series. This paper considers a popular online change-point problem of detecting a change in distribution of a sequence of i.i.d. random variables. Online change-point problems are concerned with monitoring the structure of a random process(es) whose observations arrive sequentially. For these problems, any good monitoring procedure should reliably alert the user to unexpected changes as soon as possible or with highest probability, subject to a tolerance on false alarms.

Let $y_1,y_2,\ldots$ be a sequence of independent random variables arriving sequentially. The purpose of this paper is to discuss tests for the hypothesis that $y_i \,(i=1,2\ldots)$ are identically distributed with some probability density function (pdf) $f(y)$ against the alternative that at some unknown change point $0\leq \nu < \infty$, the random variables $y_1,y_2,\ldots,y_\nu $ and  $y_{\nu+ l+1},y_{\nu+l+2},\ldots $ are identically distributed with density $f(y)$ and $y_{\nu+1},y_{\nu+2},\ldots,y_{\nu+l} $ are identically distributed with pdf $g(y)$ such that $g(y)\neq f(y)$. Here, $l$ is length of the change-point period (signal) and can be known or unknown. Under a standard hypothesis testing framework, the null hypothesis is $ \mathbb{H}_\infty\!:\, \nu = \infty $ and hence the pdf $f(y)$ is the density of $y_i$ for all $i=1,2,\ldots $. The alternative hypothesis is $ \mathbb{H}_\nu$: $0\leq \nu< l\leq \infty$ and therefore
\bea %\label{means_of_eps}
\mathbb{H}_\nu\!: \left \{\begin{array}{ll}
                           \mbox{ \!\!\!\!$y_i$ have density $f(y)$}  & \mbox{ if $i\leq \nu$ or  $i>\nu+l$ }\\
            		\mbox{ \!\!\!\!$y_i$ have density $g(y)$} & \mbox{ if $\nu < i \leq \nu+l$ }
                         \end{array} \right.
\eea
with $i=1,2,\ldots$. Under $\mathbb{H}_\nu$, the arrival time of the signal is $\nu+1$ (it is unknown). Most classical results assume $f$ and $g$ are known completely; by this, we mean no nuisance parameters are present in the distributions. In later sections, we will briefly discuss tests designed to approach the change-point problem in the presence of nuisance parameters.

A thorough introduction to the field of online (quickest) change-point detection mainly for the case of $l=\infty$ can be found in, for example, \cite{basseville1993detection,tartakovsky2014sequential,zhigljavsky1988detection,poor2008quickest}. Some of the most popular online change-point algorithms used in practice are Shewhart's $\bar{X}$-chart \cite{shewhart1931economic}, the CUSUM algorithm \cite{page1954continuous}, the Shiryaev-Roberts procedure \cite{shiryaev1963optimum,roberts1966comparison} and the Exponentially Weighted Moving Average (EWMA) chart \cite{roberts2000control}. The case of $l=\infty$, and hence when a change in distribution occurs it does so permanently, is by far the most popular scenario considered in the change-point literature; a number of influential papers are \cite{lorden1971procedures, moustakides1986optimal,pollak1985optimal, ritov1990decision,gordon1994efficient,lai1995sequential}. The CUSUM and Shiryaev Roberts procedures benefit with their simplicity and proven optimality under suitable optimality criteria; these two procedures will be the focus of discussion for the case $l=\infty$.
 The case of finite $l$, and hence when a change occurs it does so temporarily, has seen considerable attention in the past, see \cite{guepie2017detecting,guepie2012sequential,lai1974control,Chu,eiauer1978use,Bau2,glaz2009scan,Glaz2012}. More recently it has been the focus of attention in the papers of \cite{noonan2020power,noonan2020approximations,zhigljavsky2021first}.  Examples of areas where detecting a transient change in distributions is extremely important can be found in radar and sonar \cite{bell1993information,poor2013introduction,streit1999detection}, nondestructive testing \cite{schmerr2016fundamentals}, and medicine \cite{bianchi1993time}. Non-parametric online change-point detection methods have also become very popular \cite{MZ2003,brodsky2013nonparametric}.  For the state of the art techniques for multiple change-point detection, see \cite{fryzlewicz2014wild,cho2015multiple,korkas2017multiple,fryzlewicz2020detecting}. For sequential change-point detection in high-dimensional time series, a likelihood ratio approach can be found in \cite{dette2020likelihood,gosmann2019new}.

This survey is organised as follows.
In Section~\ref{sec:permanent}, we survey results for $\l=\infty$ and discuss known optimality results for the CUSUM and Shiryaev-Roberts procedures. This section contains well known classical results but is included to introduce the reader to change-point concepts that will be used when considering the transient change-point problem. In Section~\ref{sec:temporary}, we assume $l<\infty$ and discuss a number of online tests for transient changes; the likelihood ratio test providing the inspiration behind all tests. In this Section, we compare procedures when applied for detecting a temporary change in mean of a sequence of Gaussian random variables. We also apply tests for monitoring stability of components used in the Oil and Gas industry.

Throughout this survey we shall use the notation ${\rm Pr}_{\infty}$ and $\mathbb{E}_\infty$ to denote probability and expectation under $\mathbb{H}_\infty$. Under the alternative $\mathbb{H}_\nu$, we shall use the notation ${\rm Pr}_{\nu}$ and $\mathbb{E}_\nu$ to denote probability and expectation assuming the change-point occurs at $\nu<\infty$.

\section{Permanent change in distributions}\label{sec:permanent}
 In this section, we assume $l=\infty$; if a change occurs, it does so permanently. Suppose $y_1,y_2,\ldots,y_n$ have been sampled. The likelihood ratio for testing $\mathbb{H}_\infty$ against $\mathbb{H}_\nu$ is
\bea
\Lambda_{\nu,n} = \prod_{i=\nu+1}^{n}\frac{g(y_i)}{f(y_i)}
\eea
assuming $\nu<n$, otherwise $\Lambda_{\nu,n}=1$. 

\subsection{The CUSUM and Shiryaev-Roberts procedures}
By maximising the statistic $\Lambda_{\nu,n}$ over all possible locations of $\nu$, we obtain the CUSUM statistic
\be\label{CUSUM_stat}
&& V_n := \max_{1\le \nu \le n} \Lambda_{\nu,n}, \,\,\,\,\, n\geq 1 \,.
%\prod_{j=k}^{n} \exp\left( \frac{(y_j-\mu)^2-(y_j-\mu-A)^2}{2\sigma^2}  \right)\,.
\ee
The CUSUM stopping rule (when to alert the user to a potential change-point) is
\be\label{CUSUM_test}
&&\tau_{V }(H) := \inf \{n\ge 1: V_n > H  \}\,.
\ee
An appealing property of statistic \eqref{CUSUM_stat} is the recursive property
\bea
&&V_n = \max\{ V_{n-1},1 \}\cdot \frac{g(y_n)}{f(y_n)},  \,\,\,\,\, V_0= 1 \,.
\eea
The threshold $H$ in $\tau_{V }(H) $ is chosen on the users tolerance to false alarm risk. Page \cite{page1954continuous} and Lorden \cite{lorden1971procedures} measured false alarm risk through the Average Run Length to false alarm (ARL). This corresponds to choosing $H$ such that $\mathbb{E}_\infty\tau_{V}(H)=C$,  where $C$ is a pre-defined value chosen by the user but is typically large. How to compute $\mathbb{E}_\infty\tau_{V}(H)$ will be discussed later in this section.

The famous CUSUM chart of Page \cite{page1954continuous} introduces a reflective barrier at zero:
\be\label{Page_procedure}
P_n \!\!\!\!&=&\!\!\!\! \max \left\{ P_{n-1}+\log \frac{g(y_n)}{f(y_n)},0 \right\},\\ P_0\!\!\!\!&=&\!\!\!\! 0 \,. \nonumber
\ee

The statistics \eqref{Page_procedure} and $\log V_n $ are equivalent on the positive half plane and hence the  rule
\be\label{Page_CUSUM_test}
&&\tau_{P}(\log(H)) = \inf \{n\ge 1: P_n > \log H  \}\,,
\ee
and $\tau_{V}(H)$ are equivalent for $H>1$. The stopping rule $\tau_{V}$ is more general than $\tau_{P}$ as thresholds $H \leq 1$ are permissable. An approximation for $\mathbb{E}_\infty\tau_{P}(H)$ for general distributions $f$ and $g$ was derived in \cite{polunchenko2012nearly}. Let $I_f:= -\mathbb{E}_\infty (\log(g(y_1)/f(y_1)))$ and $I_g =\mathbb{E}_0(\log(g(y_1)/f(y_1))) $ (to compute $I_g$ we assume the change-point occurs at time zero). Then
\be\label{CUSUM_ARL}
&&\mathbb{E}_\infty\tau_{P}(H) \simeq \frac{e^H}{I_g\zeta^2}- \frac{H}{I_f} - \frac{1}{I_g\zeta} \,.
\ee
Here the constant $\zeta$ is called the limiting exponential overshoot. Let $Z_n = \sum_{i=1}^{n}\log(g(y_i)/f(y_i))$ be a random walk. Then it can be shown, see \cite[Ch. VIII]{Sieg_book}, that
\bea
\zeta \!=\! \frac{1}{I_g}\exp \left\{\!-\!\sum_{k=1}^\infty\frac{1}{k}[{\rm Pr}_{\infty}(Z_k\!>\!0)\!+\!{\rm Pr}_{0}(Z_k\!\leq\!0)  ]  \right\}\! .
\eea 
The approximation \eqref{CUSUM_ARL} seems extremely accurate. For example, suppose pre-change observations are i.i.d $N(0,1)$ random variables and post-change observations are i.i.d $N(A,1)$ for some known $A>0$. We have
\be\label{Gaussian_problem}
f(y) \!\!\!\!&=&\!\!\!\! \frac{1}{\sqrt{2\pi}}\exp(-y^2/2),\\
g(y) \!\!\!\!&=&\!\!\!\! \frac{1}{\sqrt{2\pi}}\exp(-(y-A)^2/2) \,. \nonumber
\ee
For $A=1$, Monte Carlo simulations provide $\mathbb{E}_\infty\tau_{P}(4.39) =500$. Application of the approximation in \eqref{CUSUM_ARL} provides 498. The draw back of the approximation in \eqref{CUSUM_ARL} is that $\zeta$ requires expensive numerical evaluation.

To construct the Shiryaev-Roberts (SR) procedure, define the generalised Bayesian detection statistic as:
\be\label{SR_stat}
&& R_n := \sum_{\nu=1}^{n}\prod_{j=\nu+1}^{n} \Lambda_{\nu,n} \,.
\ee

Then the SR test is:
\be\label{Shiryaev Roberts_test}
&&\tau_{R}(H) := \inf \{n\ge 1: R_n > H  \}\,,
\ee
where $H$ is the solution of  $\mathbb{E}_\infty\tau_{R}(H)=C$ for some pre-determined $C$. The SR statistic \eqref{SR_stat} satisfies the following recurrence:
\bea
&&R_n = (1+R_{n-1}) \cdot \frac{g(y_n)}{f(y_n)},  \,\,\,\,\, n\geq 1, \, R_0=0 \,.
\eea

\subsection{Evaluating ARL for CUSUM and SR tests}

Explicit expressions for $\mathbb{E}_\infty \tau_{V}(H)$ and  $\mathbb{E}_\infty \tau_{R}(H)$ are not known. However, they can be numerically obtained by numerically solving particular Fredholm integral equations as proved in \cite{moustakides2009numerical}. Here it was shown that $\mathbb{E}_\infty \tau_{ V}(H)$ and $\mathbb{E}_\infty \tau_{R}(H)$ can be computed by a unified approach for general Markov statistics. Set $H>0$. For a sufficiently smooth positive valued function $\xi$ and $s\in[0,H]$, let
\bea
S_{n} = \xi(S_{n-1})\cdot \frac{g(y_n)}{f(y_n)}\,\,\,\,\, n\geq 1, \,\, S_0=s \in[0,H]
\eea
be a Markov detection statistic with stopping rule
\bea
\tau_S(H):=\inf \{n\ge 1: S_n > H  \}\,.
\eea
Let $\phi(s) =\mathbb{E}_\infty(\tau_S(H)) $ be the ARL (note the dependence on $S_0=s$) and  set $F(x) = {\rm Pr}_{\infty}({g(y_1)}/f(y_1)\leq x)$. Then $\phi(s)$ is the solution of the following Fredholm integral equation:
\be\label{ARL_fredholm}
&&\phi(s) = 1 + \int_{0}^{H}\phi(x)\left[\frac{d}{dx}F\left(\frac{x}{\xi(s)} \right) \right]dx\,.
\ee
For the CUSUM  and SR procedures we have $\xi(s) = \max(1,s)$ and $\xi(s)=1+s$, respectively.  To solve this integral equation, we refer to \cite{moustakides2009numerical}.

Approximations for ARL of the CUSUM and SR procedures have been specifically developed for the problem of detecting the change in mean of normal random variables. Here we operate under \eqref{Gaussian_problem}.  
To approximate ARL for both the CUSUM and SR procedures or to narrow the domain of search and more efficiently numerically solve the Fredholm equation \eqref{ARL_fredholm}, one could use the following simple approximations developed in \cite{tartakovsky2005asymptotic} and \cite{pollak1987average} respectively:
\be
\mathbb{E}_\infty\tau_{V}(H_{}) &\simeq& 2H_{}/(A\kappa^2(A)) \, , \label{tart_app}\\
\mathbb{E}_\infty\tau_{R}(H_{}) &\simeq&H_{}/\kappa(A)  \label{poll_app}\,,
\ee
where
\bea
\kappa(A) = \frac{2}{A^2}\exp \left\{ -2\sum_{\nu=1}^{\infty}\frac{1}{\nu}\Phi\left(-\frac{A}{2}\sqrt{\nu} \right) \right\} \\\mbox{ and }\Phi(x) = \int_{-\infty}^{x}f(y)dy \,.
\eea

The approximations in \eqref{tart_app} and \eqref{poll_app} are extremely accurate. In Table~\ref{tart_table}, one can observe the high accuracy of approximation  \eqref{tart_app} for different thresholds $H_{}$. In fact, \eqref{poll_app} is remarkably accurate and frequently leads to exact values of ARL. The only slight inconvenience of both approximations is the numerical evaluation required to compute $\kappa(A)$. This quantity can be approximated with $\kappa(A) \simeq \exp(-\rho\cdot A)$, where the constant $\rho$ is defined later in \eqref{D_rho} but can be approximated to three decimal places by $\rho\simeq0.583$. Using this approximation for $\kappa$ in \eqref{tart_app} and  \eqref{poll_app} still results in excellent approximations.

\begin{table*}
\caption{Approximations for $\mathbb{E}_\infty\tau_{V}(H)$  with $A=1$.}
\label{tart_table}
\centering
\begin{tabular*}{\textwidth}{@{\extracolsep{\fill}}cccccc}
\hline
$H$ & 9.32& 17.33 & 80.65 &159.35 & 788.00   \\
\hline
$\mathbb{E}_\infty\tau_{V}(H)$ & 50 & 100 & 500 &1000   & 5000            \\
 Approximation \eqref{tart_app}& 59 & 110 & 513  &1014  & 5018      \\
 Approximation \eqref{tart_app} with $\kappa(A) \simeq \exp(-\rho\cdot A)$& 60 & 111 & 517  &1023  & 5058      \\
\hline
\end{tabular*}
\end{table*}

%The power of the CUSUM test for the transient change considered in this paper is defined as
%\bea
%\!\!\!\!\! {\cal{P}}_{CS}(H_{CS},A)\! &:=&\!\!\! \lim_{\nu \rightarrow \infty}{\rm P}_{1}\{V_{n}\!>\!H_{CS}  \text{  for some }  n\!\in\![\nu+\!1,\nu\!+\!2l\!-\!1  ]  \, | \, \tau_{CS}(H_{CS}) \!>\! \nu  \}.
%\eea
%To compare the MOSUM test with the CUSUM test, the threshold $H$ in \eqref{def:MOSUM} and the threshold $H_{CS}$ in \eqref{CUSUM_test} are chosen such that $\mathbb{E}_0\tau(H)=\mathbb{E}_0\tau_{CS}(H_{CS})=5000$. Determination of $H_{CS}$ for CUSUM was obtained using tabulated values given in \citet[p. 3237]{moustakides2009numerical} whereas  suitable values of $H$ for MOSUM was obtained using  Approximation 0 of Section~\ref{ARL_section}.
%
%

\subsection{Optimality criteria} \label{Comparing_procedures}

Denote by $\Delta(C)$ the set of all stopping times of change-point procedures with ARL of at least $C$. More precisely,
$
\Delta(C) := \{\tau: E_{\infty}\tau \geq C \},\,\,\, C>1,
$
where $\tau=\tau(H)$ is a stopping time for a sequential change-point procedure.
A common criterion for comparing change-point procedures when $l=\infty$ is the supremum Average Delay to Detection (ADD) introduced by Pollak \cite{pollak1985optimal}. Define $ADD_\nu(\tau) := \mathbb{E}_\nu(\tau-\nu|\tau>\nu)$. Then
\be\label{Pollak_sad}
&&SADD(\tau) := \sup_{0\leq \nu<\infty}ADD_\nu(T) \,.
\ee
An optimal change-point procedure would satisfy $SADD(\tau_{opt}) = \inf_{\tau \in \Delta(C)}SADD(\tau)$ for all $C>~1$. Finding an optimal procedure for this criterion is very difficult, where in general only asymptotic optimality as $C\rightarrow \infty$ (low false alarm rate) is known \cite{pollak1985optimal}. Another popular criterion is the worst-case minimax scenario of Lorden \cite{lorden1971procedures}  defined as 
\be\label{lorden}
{\cal L}(\tau):=\sup_{\nu\geq 0}\esssup \mathbb{E}_\nu[ (\tau-\nu)^+|y_1,y_2,\ldots,y_\nu ] \,. \nonumber\\
\ee
This criterion evaluates the average detection delay conditioned on the worst possible data before the change and then considers the worst possible deterministic change-point. Asymptotic optimality (as $C\rightarrow \infty$) of the CUSUM chart of Page was proved in \cite{lorden1971procedures}. It was subsequently proved in \cite{moustakides1986optimal} that the CUSUM chart of Page is in fact optimal under this criterion for every $C>1$.

The SR procedure is optimal for every $C>1$ under the Stationary Average Delay to Detection (STADD) criterion. The STADD criterion praises detection procedures that detect the change as quickly as possible, at the expense of raising many false alarms (using a repeated application of the same stopping rule). Formally, the STADD criterion is defined as follows. Let $\tau_1,\tau_2\ldots$ be a sequence of independent copies of the stopping time $\tau$. Let $T_j=\tau_1+\tau_2+\ldots+\tau_j$ be the time the $j^{th}$ alarm is raised. Let $I_\nu = \min\{j>1:T_j>\nu\}$; this is the index of the first alarm which is not false after $I_\nu-1$ false alarms. Then
\bea
&&STADD(\tau) := \lim_{\nu \rightarrow \infty} \mathbb{E}_\nu[T_{I_{\nu}}-\nu] \,.
\eea

The STADD criterion is equivalent to the Relative Integral Average Detection Delay (RIADD) measure, see \cite{moustakides2009numerical}, which is defined as:
\bea
&&RIADD(\tau) = \frac{\sum_{\nu=0}^{\infty}\mathbb{E}_\nu[(\tau-\nu)^+ ]}{\mathbb{E}_\infty[\tau]} \,.
\eea

Both the SR procedure and CUSUM procedure are asymptotically optimal as $C\rightarrow \infty$ for the Lorden and the STADD criteria. It is discussed in \cite{moustakides2009numerical} for both CUSUM and the Shiryaev–Roberts procedure Lorden’s essential supremum measure \eqref{lorden} and Pollak’s supremum measure $SADD$ defined in \eqref{Pollak_sad} are attained at $\nu=0$, that is:
\bea
&&{\cal L}(\tau_{V}(H))=SADD(\tau_{V}(H)) = \mathbb{E}_0\tau_{V}(H), \\
 &&{\cal L}(\tau_{R}(H))=SADD(\tau_{R}(H)) = \mathbb{E}_0\tau_{R}(H) \,.
\eea

Similarly to the computation of $\mathbb{E}_\infty\tau_{V}(H)$ and $\mathbb{E}_\infty\tau_{R}(H)$, to obtain $\mathbb{E}_0\tau_{V}(H)$ and $\mathbb{E}_0\tau_{R}(H)$ one can numerically solve a Fredholm equation.  Instead of setting
 $\phi(s) =\mathbb{E}_\infty(\tau(H)) $, let $\phi(s) =\mathbb{E}_0(\tau(H)) $. Also set $F(x) = {\rm Pr}_{0}({g(y_1)}/f(y_1)\leq x)$. Then from \cite{moustakides2009numerical}, $\phi(s)$ is the solution of the Fredholm integral equation given in \eqref{ARL_fredholm}. The computation of STADD requires solving a slightly more difficult integral equation and we refer the interested reader to \cite{moustakides2009numerical} for more discussions.
 
For the Gaussian example considered in \eqref{Gaussian_problem}, the findings of \cite{moustakides2009numerical} indicate that for small values of $A$ say $A=0.01$, the CUSUM noticeably outperforms the SR procedure under Lordens criterion. Vice versa, the SR procedure noticeably outperforms CUSUM under the STADD framework. When the change in $A$ becomes large, say $A=1$, the benefits a procedure has over the other diminishes.

\section{Transient change in distributions}\label{sec:temporary}
In this section, we assume $1\leq l<\infty$ and therefore study procedures aimed at detecting a transient change in distributions. Suppose $y_1,y_2,\ldots,y_n$ have been sampled.
The log likelihood ratio for testing $\mathbb{H}_\infty$ against $\mathbb{H}_\nu$  is
\be\label{MOSUM_likelihood}
&&\Lambda_{\nu,\nu+l} = \sum_{i=\nu+1}^{\min\{\nu+l,n\}} \log \frac{g(y_i)}{f(y_i)} \,.
\ee

\subsection{A collection of procedures}

%The Levin and Kline statistic \cite{levin1985cusum} is

For $l$ unknown, the log likelihood ratio statistic is obtained by maximising \eqref{MOSUM_likelihood} over all possible change point locations $\nu$ and transient change lengths:
\be\label{Levin_statistic}
&&K_n := \max_{0\leq \nu <\nu+l \leq n}\Lambda_{\nu,\nu+l}\,,
\ee
with the stopping rule 
\be\label{LK_stopping}
&&\tau_{K}(H):=\inf \{n\ge 1: K_n > H  \}\,.
\ee
%If nuisance p
If there are no nuisance parameters present in $f$ and $g$ that require estimation, the statistic \eqref{Levin_statistic} satisfies the recursive property:
\be\label{likelihood_recursion}
K_n \!\!\!\!&=&\!\!\!\! \max\{ K_{n-1},\max_{0\leq\nu\leq n-1}\Lambda_{\nu,n} \},  \\
 K_0\!\!\!\!&=&\!\!\!\! 0 \,. \nonumber
\ee

For large $n$, the statistic \eqref{Levin_statistic} is very expensive to compute despite the recursive property given in \eqref{likelihood_recursion}. For offline change-point problems, this large computational expense may be an inconvenience but it is not a fundamental problem as time is often not an issue. However, for online procedures that require calculations in real time, the statistic $K_n$ is not practical. The assumption that no knowledge of the transient change length is known is unlikely. One can imagine that some prior knowledge about the length of transient change is likely, for example it may be bounded $l_0\leq l \leq l_1$. From here on, this assumption will be made. The log likelihood ratio statistic is:
\be\label{general_LK}
&&Z_n={Z}_n(l_0,l_1) := \max_{\substack{0\leq \nu <\nu+l \leq n\\l_0 \leq l \leq l_1} }  \Lambda_{\nu,\nu+l}\,,
\ee
with the stopping rule 
\bea
\tau_{{Z}}(H):=\inf \{n\ge l_1: {Z}_n(l_0,l_1) > H  \}\,.
\eea
If no nuisance parameters require estimation, the statistic ${Z}_n$ satisfies the following recursive property:
\bea
{Z}_n \!\!\!\!&=&\!\!\!\! \max\{ {Z}_{n-1},\max_{n\leq\nu\leq n+l_1-l_0}\Lambda_{\nu,n+l_1} \} , \\
 Z_{l_1} \!\!\!\!&=&\!\!\!\! \max_{\substack{0\leq \nu <\nu+l \leq l_1\\l_0 \leq l \leq l_1} }  \Lambda_{\nu,\nu+l} \,.
\eea
This is much easier to compute than \eqref{Levin_statistic} for  $n$ large. This recursive property means the stopping rule of $Z_n$ can be expressed as:
\be\label{stopping_rule_bounded}
\tau_{Z}(H)\!\!\!\!&=&\!\!\!\! l_1+\tau_{{{S,l_0,l_1}}}(H)\,, \text{ where }\\
\tau_{{{S,l_0,l_1}}}(H)\!\!\!\!&:=&\!\!\!\! \inf \{n\ge0: {S}_{n,{l_0,l_1}} > H  \}\,, \nonumber
\ee 
and
\bea
{S}_{n,{l_0,l_1}} :=\max_{n\leq\nu\leq n-l_0+l_1} \Lambda_{\nu,n+l_1} \,.
\eea
Therefore  $\mathbb{E}_\infty\tau_{Z}(H)=\mathbb{E}_\infty\tau_{{{S,l_0,l_1}}}(H)+l_1$\,.\\

If we make the additional assumption that $l$ is known exactly and is completely contained within the sample of size $n$, i.e. $\nu+l\leq n$, then the MOSUM statistic can be obtained by setting $l_0=l_1=l$ in \eqref{general_LK}. For this reason, the statistic $Z_n$ can be called the generalised MOSUM procedure. This corresponds to maximising \eqref{MOSUM_likelihood} over all valid change-point locations $\nu$:
\be\label{Levin_kline_stat}
&&M_n := \max_{0\leq \nu \leq n-l}\Lambda_{\nu,\nu+l} \,,
\ee
with the stopping rule
\bea
&&\tau_{M}(H):=\inf \{n\ge l: M_n > H  \}\,.
\eea

In what follows, we will define the MOSUM test for a general window length $L$, with $L$ a fixed positive integer. Results for the likelihood ratio test can be obtained by setting $L=l$. Define the moving sums
\begin{eqnarray*}
S_{n,L}:={S}_{n,{L,L}}= \sum_{j=n+1}^{n+L} \log \frac{g(y_i)}{f(y_i)}\, \;\; (n=0,1, \ldots)\,. %\label{eq:sumsq2}.\,.
\end{eqnarray*}
Then the stopping rule $\tau_{M}(H)$ for a given window length $L$ can be expressed as
\be\label{MOSUM_stopping}
\tau_{M}(H) \!\!\!\!&=&\!\!\!\! \tau_{S,L}(H)+L,\, \mbox{ where }\\
 \tau_{S,L}(H)\!\!\!\!&:=&\!\!\!\!\inf \{n\ge 0:  S_{n,L} > H  \}\, \nonumber
\ee
and therefore $\mathbb{E}_\infty\tau_{M}(H)=\mathbb{E}_\infty\tau_{S,L}(H)+L$. The moving sum $S_{n,L}$ is the reason behind the MOSUM name.\\

%By definition, the ARL $\mathbb{E}_\infty \tau_{\tilde{K}}(H)$ satisfies:
%\bea
%\mathbb{E}_\infty {\tau_{\tilde{K}}}(H)&=& l_1+\sum_{k=1}^\infty k \cdot {\rm P}_\infty \{{\tau_{\tilde{S}}}(H)=k\} \,\nonumber \\
%&=&l_1+\sum_{k=1}^\infty k \left[ {\rm P}_\infty \{{\tau_{\tilde{S}}}(H)\leq k\}-{\rm P}_\infty \{{\tau_{\tilde{S}}}(H)\leq k-1\}  \right] \, \nonumber \\
%&=& l_1+\sum_{k=1}^\infty k \left[  {\rm P}_\infty\left\{ \max_{j=0\ldots k}{S}_{j,{l_0,l_1}}\geq H \right\}-  {\rm P}_\infty\left \{ \max_{ j=0\ldots k-1}{S}_{j,{l_0,l_1}}\geq H\right\}  \right] \label{ARL}
%\eea
% 
%
%Therefore, computation or approximations for  $\mathbb{E}_\infty \tau_{\tilde{K}}(H)$ can be obtained by approximating the probability ${\rm P}_\infty\left\{ \max_{j=0\ldots k}\tilde{S}_{j,{l_0,l_1}}\geq H \right\}$ for all $k$ and $H$. By similar arguments, $\mathbb{E}_\infty {\tau_{M}}(H)$ can be expressed as:
%\be\label{MOSUM_ARL_Def}
%\mathbb{E}_\infty {\tau_{M}}(H)= l + \sum_{k=1}^\infty k \left[ {\rm P}_\infty\left(\max_{i=0\ldots,k}S_{i,l}\geq H \right)-  {\rm P}_\infty\left(\max_{i=0\ldots,k-1}S_{i,l}\geq H \right) \right] \,.
%\ee

For the transient change-point problem, the false alarm risk can be measured through ARL. However, this is not the only approach taken in the change-point literature. In \cite{guepie2012sequential} and \cite{lai1998information} , the false alarm risk is measured through:
\be\label{igor_false_alarm}
&&\sup_{k\geq 1}{\rm Pr}_\infty(k\leq \tau<k+m_\alpha)\leq \alpha \,,
\ee
where $\tau$ is a stopping rule, $\alpha$ is your false alarm tolerance (type 1 error) and $\lim \inf m_\alpha/|\log(\alpha)|>I_g^{-1}$ but $\log m_\alpha=o(| \log\alpha| )$ as $\alpha \rightarrow 0$; recall $I_g =\mathbb{E}_0(\log(g(y_1)/f(y_1))) $. Another alternative to the usual ARL constraint has been proposed in \cite{tartakovsky2005asymptotic,tartakovsky2008discussion}. Here, the suggested criterion is
\bea
&&\sup_{k\geq 1 }{\rm Pr}_\infty(\tau < k+ m_\alpha | \tau \geq l )\leq \alpha.
\eea

From now on, false alarm risk will be measured through ARL $\mathbb{E}_\infty \tau$ and we refer to \cite{guepie2012sequential,lai1998information,tartakovsky2005asymptotic,tartakovsky2008discussion} for more discussions on other approaches. The majority of research has focused on detecting transient changes in a sequence of Gaussian random variables. The next section is devoted solely this problem. 
\subsection{Detecting a transient change in Gaussian random variables}

Consider the problem of detecting the change in mean of normal random variables. Suppose pre-change observations are i.i.d $N(\mu,1)$ random variables and post-change observations are i.i.d $N(\mu+A,1)$ for some $A>0$. The values of $\mu, l$ and $A$ may be known or unknown, with $\mu$ and $A$ playing the roles of nuisance parameters if unknown.  We have
\be\label{Gaussian_problem2}
f(y) \!\!\!\!&=&\!\!\!\! \frac{1}{\sqrt{2\pi}}\exp(-(y-\mu)^2/2), \\
g(y) \!\!\!\!&=&\!\!\!\! \frac{1}{\sqrt{2\pi}}\exp(-(y-\mu-A)^2/2) \,. \nonumber
\ee

The offline version of this change-point problem is devoted to testing for change-points in a sample of fixed length and has seen significant attention in the past, see \cite{levin1985cusum,siegmund1986boundary,hogan1986large,siegmund1988approximate,yao1989large}. An excellent survey of several statistics aimed at addressing the offline problem can be found in \cite{yao1993tests}. Despite the fact $Z_n$ defined in \eqref{general_LK} is a generalisation of $M_n$ given in \eqref{Levin_kline_stat}, we will initially discuss recent results for $M_n$. These results will provide inspiration for addressing the much more complicated problems associated with $Z_n$.

\subsubsection{The MOSUM statistic}\label{sec:Mosum}
%Here we make the assumption $l$ and $\mu$ are known in \eqref{Gaussian_problem}. The case of $l,\mu$ and $A$ unknown can be obtained in the next section as a special case by setting $l_0=l_1$. 
For the MOSUM test, knowledge of $A$ is not required to set the ARL constraint; this is because the MOSUM stopping rule given in \eqref{MOSUM_stopping} specialised for this Gaussian example is tantamount to:
\be\label{MOSUM_normal}
\tau_{M}(H)\!\!\!\!&=&\!\!\!\! \tau_{S,L}(H)+L, \\
 \tau_{S,L}(H)\!\!\!\!&=&\!\!\!\!\inf \{n\ge 0: S_{n,L} > H  \}\, ,  \nonumber \\
\text{ with } S_{n,L}\!\!\!\!&=&\!\!\!\! \sum_{j=n+1}^{n+L} y_j\,  \,. \nonumber
\ee

%\subsection{Approximation for $\mathbb{E}_0 {\tau}_{\xi}(h$) }
The problem of approximating  $\mathbb{E}_\infty \tau_{S,L}(H)$ assuming $\mu$ is known was considered {in \cite{noonan2020approximations}}. Here will recall the main steps in the construction.
Define
 \bea
\label{H_h}
h= \frac{H - \mu L}{ \sqrt{L}}
\;\mbox{ so that }\;H= \mu L +
h\sqrt{L} \, 
\eea 
and consider the standardised versions of $S_{n,L}$:
\bea
\label{eq:def-xi}
\xi_{n,L}\!:=\! \frac {S_{n,L}-  \E_\infty  S_{n,L}}
{\sqrt{{\rm Var}_\infty(S_{n,L})}}\!=\!
\frac { S_{n,L}- \mu L}
{ \sqrt{  L    }}  \, ,\mbox{$n=0,1,\ldots\, .$}
\eea

Then the stopping time $\tau_{S,L}(H)$ is equivalent to the stopping time
\be\label{relation2}
&&{\tau}_{\xi}(h) := \inf \{ n\ge 0: \xi_{n,L} \ge h \}
\ee
 and hence $\mathbb{E}_\infty{\tau}_{\xi}(h) =\mathbb{E}_\infty\tau_{S,L}(H)$.

For any integer $M\ge0$, the discrete time process $\xi_{0,L},\xi_{1,L},\ldots, \xi_{M,L}$  is approximated by a continuous time analogue $S(t)$ on $[0,T=M/L]$. The process $S(t)$  is a zero mean, stationary Gaussian process with correlation function $R(t) = \max \{ 0,1-|t| \}$. The ARL  $\mathbb{E}_\infty{\tau}_{\xi}(h)$ then has the continuous-time approximation
\be\label{ARL_cont_approx}
&&\mathbb{E}_\infty {\tau}_{\xi}(h)\cong-L \int_0^{\infty}{s \,dF_h(s)}\, ,
\ee
where  $F_{h}(T) := {\rm P}_\infty(S(t)< h \text{  for all   } t\in[0,T])$.

Explicit formulas for the probability $F_{h}(T)$ with $T\le 1$ were first derived in \cite{slepian1961first}. Here it was shown
\bea
F_h(T)\!\!\!\!&=&\!\!\!\!\int_{-\infty}^{h}\!\!\Phi\left(
\frac{h({ Z}\!+\!1)\!-\!x(-{ Z}\!+\!1)}{2\sqrt{{ Z}}}\right)\varphi(x)dx \nonumber\\
&&\!\!\!\!\!\!\!\!\!\!\!\!\!\!\!\!\!\!\!\!\!\!\!\!\!\!\!\!- \frac{2\sqrt{Z}}{Z+1} \varphi(h)\left[{h {\sqrt{Z}}}\,\Phi(h\sqrt{{ Z}})
\!+\! \frac1{ \sqrt{2\pi} } (\sqrt{2\pi} \varphi(h))^{Z} \right]\,  \nonumber.
\eea
For $T=1$ this reduces to
\be \label{one_step}
&&F_h(1)=\Phi^2(h)-\varphi(h) \big[h\Phi(h)+ \varphi(h) \big] .
\ee

 For $T>1$, formulae for $F_{h}(T)$ were first derived in \cite{Shepp71}; these expressions take different forms depending on whether or not $T$ is integer.
The result of  \cite[p.949]{Shepp71} states than if $T=n$ is a positive integer then
\be \label{Shepp_form}
{ F}_{h}(n) = \int_{-\infty}^{h}\int_{D_x} \det[\varphi(y_i - y_{j+1} + h)]^n_{i,j=0} \,\nonumber  \\
dy_2\ldots dy_{n+1}dx\,
\ee where $y_0= 0, y_1=h-x,$
$
D_x =  \{y_2, \dots , y_{n+1} \given h-x < y_2 < y_3 < \ldots < y_{n+1}   \}
$.
For non-integer $T\ge1$,  the exact formula for $
{ F}_{h}(T)
$ is even more complex (the integral has the dimension $\lceil 2T \rceil +1$)  see \cite[p.950]{Shepp71}.
For $T=2$, \eqref{Shepp_form} yields
\be\label{Prob_2}
F_h(2) \!\!\!\!&=&\!\!\!\! \Phi^3(h) - 2h\varphi(h) \Phi^2(h) \nonumber   \\
\!\!\!\!&+& \!\!\!\!\frac{h^2\!-\!3\!+\!\sqrt{\pi} h}{2}   \varphi^2(h)\Phi(h)+ \frac{h \!+\! \sqrt{\pi}}{2}\varphi^3(h) \nonumber\\
\!\!\!\! &+& \!\!\!\!\int_{0}^{\infty}\Phi(h-y) [\varphi(h+y)\Phi(h-y) \nonumber \\
\!&& - \sqrt{\pi}\varphi^2(h)\Phi(\sqrt{2}y) \, ] dy.
\ee

The complicated nature of these expressions for $F_{h}(T)$ made them impractical for the use in the ARL approximation \eqref{ARL_cont_approx}. One simple yet still very accurate approximation has the form (see \cite{noonan2018approximating}):
\be\label{F_approx}
&&F_h(T) \simeq F_h(2) \left[ \theta(h)  \right]^{T-2} \,, 
\ee
where  $\theta(h) = {F_h(2)}/{ F_{h}(1)}$ and the probabilities $F_h(1  )$ and $F_h(2 ) $ are given in \eqref{one_step} and \eqref{Prob_2} respectively.
Here, $\varphi(x)$ and $\Phi(x)$ are the standard normal density and distribution functions respectively. The approximation given in \eqref{F_approx}  applied to \eqref{ARL_cont_approx} results in the following continuous-time ARL approximation:
\bea %\label{ARL_approx_explicit}
\mathbb{E}_\infty {\tau}_{\xi}(h)  &\simeq&  -\frac{L\cdot   { F}_{h}(2)}{\theta(h)^2\log( \theta(h))}.
\eea
This approximation was then corrected in \cite[Section 7]{noonan2020approximations} for discrete time to improve results for small $L$. This amounted to correcting the probabilities ${F}_{h}(1)$ and ${F}_{h}(2)$ for discrete time; this  was performed by specialising results of D. Siegmund; primarily on expected overshoot a discrete time normal random walk has over a threshold. From \cite[p. 225]{Sieg_book}, this expected overshoot was computed as 
\be
\label{D_rho}
\rho \!\!\!\!&:=&\!\!\!\!  - \int_{0}^{\infty} \frac{1}{\pi \lambda^2}  \log\{2(1-\exp(-\lambda^2/2))/ \lambda^2 \} \, d\lambda\, \nonumber  \\
\!\!\!\!& \simeq&\!\!\!\! 0.582597. \nonumber \\
\ee
Define the probability
\be
&&{F  }_{h}(M; L) := {\rm Pr}\left(\max_{n=0,1,\ldots,M} \xi_{n,L} < h\right).
 \label{eq:prob-S}
\ee
{\ From \cite[p. 18]{noonan2020approximations}: 
\be\label{MOSUM_ARL}
\mathbb{E}_\infty \tau_{S,L}(H)=\mathbb{E}_\infty {\tau}_{\xi}(h)   \simeq  -\frac{L\cdot  {F  }_{h}(2L; L) }{\theta_L(h)^2\log( \theta_L(h))} \nonumber \\
\text{ with } 
\theta_L(h) = \frac{{F  }_{h}(2L; L)}{{ {F  }_{h}(L; L)}}\,, \nonumber \\
\ee
 where, for $h_L:=h+\omega_L$ with  $\omega_L = \sqrt{2}\rho /\sqrt{L}$, the probabilities ${F  }_{h}(L; L)$ and ${F  }_{h}(2L; L)$ can be approximated by:
\be\label{F_Corrected_shepp_explicit}
 && {F  }_{h}(L; L) \simeq \Phi(h)\Phi({h_L}) - \varphi({h_L})[h\Phi(h)+\varphi(h)]\, , \nonumber \\
\ee
\be \label{F_Corrected_shepp_explicit_2}
{F  }_{h}(2L; L) \!\!\!\! &\simeq&\!\!\!\!
  \frac{\varphi^2(h_L)}2 \!  [ ({h}^{2}\!-\!1\!+\!\sqrt{\pi}h)\Phi \left( h \right) \nonumber \\&+&(h\!+\!\sqrt{\pi})\varphi \left( h \right)
  ]    \nonumber \\
  &-&\!\varphi \left( h_L \right) \Phi \left(  h_L \right)  \left[  \left( h\!+\!
h_L  \right) \Phi \left( h  \right)\! +\!\varphi \left( h \right)  \right]\nonumber \\
 &+&\!\Phi \left( h \right)   \Phi^2({ h_L} ) \nonumber \\
  &+&\!\! \int_{0}^{\infty} \!\Phi(h\!-\!y)[\varphi(h_L+y)\Phi(h_L-y) \nonumber \\
  &&-\sqrt{\pi}\varphi^2(h_L)\Phi(\sqrt{2}y) \, ] dy.\,
  \ee
  }
Only a one-dimensional integral has to be numerically evaluated for approximating  ${F  }_{h}(2L; L)$. 
Tables~\ref{expected_run_length} and Tables~\ref{SD_run_length} demonstrate that \eqref{MOSUM_ARL} using \eqref{F_Corrected_shepp_explicit} and \eqref{F_Corrected_shepp_explicit_2} is extremely accurate.% The values of $\mathbb{E}_0 \tau_S(H)$ have been obtained by Monte Carlo simulations.

\begin{table*}
\caption{Approximations for $\mathbb{E}_\infty {\tau}_{\xi}(h)$  with $L=10$.}
\label{expected_run_length}
\centering
\begin{tabular*}{\textwidth}{@{\extracolsep{\fill}} c c c c c c c c }
\hline
$h$ & 2 & 2.25 &2.5 &2.75 & 3 &3.25 & 3.5  \\
\hline
 \eqref{MOSUM_ARL}& 126 &217  &395   & 759  &  1551  & 3375  &  7837           \\
$\mathbb{E}_\infty {\tau}_{\xi}(h)$  & 127 & 218 & 396  &757  & 1550  & 3344  &7721      \\
\hline
\end{tabular*}
\end{table*}

\begin{table*}[h]
\caption{Approximations for $\mathbb{E}_\infty {\tau}_{\xi}(h)$ with $L=50$.}
\label{SD_run_length}
\centering
\begin{tabular*}{\textwidth}{@{\extracolsep{\fill }} c c c c c c c c }
\hline
$h$ & 2 & 2.25 &2.5 &2.75 & 3 &3.25 & 3.5   \\
\hline
 \eqref{MOSUM_ARL}&  471 & 791  & 1392   & 2587  & 5099    &  10695  &    23918           \\
$\mathbb{E}_\infty {\tau}_{\xi}(h) $ & 472  &792  & 1397  &  2588 & 5085  & 10749   &    24131      \\
\hline
\end{tabular*}
\end{table*}

For approximating the boundary-crossing probability ${F  }_{h}(M; L)$ for all $M$, the discrete time corrected form of \eqref{F_approx} suggests using the approximation
\be\label{boundary_crossing_approx}
&& {F  }_{h}(M; L)\simeq  {F  }_{h}(2L; L) \left[ \theta_L(h)  \right]^{M/L-2} \,.
\ee
One could then approximate ${F  }_{h}(2L; L)$ and $\theta_L(h)$ using \eqref{F_Corrected_shepp_explicit} and \eqref{F_Corrected_shepp_explicit_2}; the high accuracy of the resulting approximation was comprehensively studied in \cite{noonan2020approximations}. 

%In Figures... we briefly show the accuracy.
%
%
%
%
%\begin{figure}[h]
%\begin{center}
% \includegraphics[width=0.5\textwidth]{L_5_T_1_CS.png}$\;\;$\includegraphics[width=0.5\textwidth]{L_5_T_2_CS.png}
%\end{center}
%\caption{Empirical probabilities of reaching the barrier $h$ (dashed black) and corresponding versions of Approximation 2 (solid red).
%Left: $T=1$ with (a) $L=M=5$ and (b) $L=M=100$. Right: $T=2$ with (a) $L=5$, $M=10$ and (b)  $L=100$, $M=200$ . }
%\label{L_5_CSA}
%\end{figure}
%

\subsubsection{The stopping rule  $\tau_{Z}(H)$ }
Here, we assume $l$ is not known exactly but can be bounded between $l_0$ and $l_1$. We will initially assume  $\mu$ and $A$ are known. The stopping rule given in \eqref{stopping_rule_bounded} specialised for this Gaussian example is tantamount to:
\be\label{generalised_MOSUM}
\tau_{{Z}}(H)\!\!\!\!&=&\!\!\!\!l_1+\tau_{{{S,l_0,l_1}}}(H)\,, \text{ where } \\
\tau_{{{S,l_0,l_1}}}(H)\!\!\!\!&:=&\!\!\!\! \inf \{n\ge0: {S}_{n,{l_0,l_1}} > H  \}\,,\nonumber 
\ee
and
\bea
{S}_{n,{l_0,l_1}} :=\max_{n\leq\nu\leq n-l_0+l_1} \sum_{j=\nu+1}^{n+l_1} \left( y_j-\mu-\frac{A}{2} \right) \,.
\eea

The short memory of the MOSUM statistic is paramount to the form of the approximation given in \eqref{boundary_crossing_approx}. This short memory is also present for the generalised moving sum statistic ${S}_{n,{l_0,l_1}}$ and suggests the form of approximations \eqref{MOSUM_ARL} and \eqref{boundary_crossing_approx} would also be suitable when applied to ${S}_{n,{l_0,l_1}}$. Introduce the probability:
\bea
{ F}_{l_0,l_1}(H, M):={\rm Pr}\{{S}_{j,{l_0,l_1}}< H \,\,\, \forall\,\,\,  j=0\ldots M\}\,.
\eea
Then the following approximations should also provide high accuracy:
\be\label{boundary_crossing_approx2}
&&\,\,\,\,\,\,{ F}_{l_0,l_1}(H, M)\!\simeq \! { F}_{l_0,l_1}(H, 2l_1) \left[ \theta_{l_1}(H)  \right]^{M/l_1-2}  \\ 
&&\,\,\,\,\, \text{ with } 
\theta_{l_1}(H)=\frac{{ F}_{l_0,l_1}(H, 2l_1)}{{ F}_{l_0,l_1}(H, l_1)} \,,\nonumber
\ee
\be\label{MOSUM_ARL2}
&& \mathbb{E}_\infty \tau_{S,l_0,l_1}(H)  \simeq  -\frac{l_1\cdot{ F}_{l_0,l_1}(H, 2l_1)   }{[\theta_{l_1}(H)]^2\log( \theta_{l_1}(H))} \,.
\ee
Unfortunately, the probability ${ F}_{l_0,l_1}(H; M)$  is complex and to the authors knowledge no formula or approximations are known. The probabilities ${ F}_{l_0,l_1}(H; 2l_1)$ and ${ F}_{l_0,l_1}(H; l_1)$ can be approximated via simulations; this is not too cumbersome as at most $3 l_1$ random variables need to be simulated at each iteration. As commonly $\mathbb{E}_\infty\tau_{{{S,l_0,l_1}}}(H)=C$ with $C$ large, say $C=500$, the right tail of the distribution of the random variable $\max_{ j=0\ldots M}{S}_{j,{l_0,l_1}}$ is of the most interest. Large deviation theory, see \cite{siegmund1986boundary,hogan1986large}, could be used to approximate the right tail of this distribution, however numerical results indicate approximations of these kind would not be accurate enough for general $l_0$ and $l_1$ (those that are not astronomically large). If the prior knowledge that $1\leq l\leq l_1$ is known, and an explicit formula to approximate ${ F}_{1,l_1}(H; M)$ or $\mathbb{E}_\infty \tau_S(H)$ is desired, the following simple heuristic argument could be used. 
The continuous time analogue of the probability ${ F}_{1,l_1}(H; M)$ is:
\bea
{\rm Pr} \left \{ \max_{\substack{0\leq s<t\leq M+l_1\\ 0\leq t-s \leq l_1}} \left[ W(t)-W(s)-\frac{A}{2}(t-s)\right] <H  \right \} \,, \nonumber \\
\eea
\vspace{-1cm}
\be\label{Brownian_prob}
\ee
where $W(t)$, $0\leq t<\infty$, is standard Brownian motion. Ideally, a large deviation approximation for \eqref{Brownian_prob} should be computed. However, for $M=l_1$, $M=2l_1$ and $A$ large, say $A\geq 1$, simulation studies indicate that the additional maximisation constraint in \eqref{Brownian_prob} of $0<t-s<l_1$ has very little influence on this probability. If this constraint is ignored, the following large deviation result of \cite{hogan1986large} can be applied.
\begin{lemma}
Suppose $\gamma>0$, $m\rightarrow \infty  $  and $u \rightarrow \infty$ such that $m\gamma u^{-1}$ is some fixed number in $(1,\infty)$. Then
\be\label{hogan_prob}
{\rm Pr} \left \{ \max_{0\leq s<t\leq m} \left[ W(t)-W(s)-\gamma(t-s)\right] >u  \right \}  \nonumber \\
= [ 2\gamma(m\gamma-u)+3+o(1)  ]\exp(-2\gamma u) \,. \nonumber \\
\ee
\end{lemma}
To subsequently correct this result for discrete time, it is recommended in \cite{hogan1986large} to increase the barrier by $H$ by $2\rho$, where $\rho$ is defined in \eqref{D_rho}. This results in the approximations
\be\label{discrete_hogan}
{ F}_{1,l_1}(H, l_1)\!\!\!\!&\simeq&\!\!\!\!1- (A(Al_1 -H-2\rho)+3 )\times \nonumber \\
&\times &\!\!\!\! \exp\{-A(H+2\rho)\}\\
{ F}_{1,l_1}(H, 2l_1)\!\!\!\!&\simeq&\!\!\!\! 1- (A(3Al_1/2 -H-2\rho)+3 )\times \nonumber  \\
&\times&\!\!\!\!\exp\{-A(H+2\rho)\}\,.
\ee
As a result, using the approximations given in  \eqref{boundary_crossing_approx2} and \eqref{MOSUM_ARL2}:
\be\label{Approx1}
{ F}_{1,l_1}(H, M)\!\!\!\!&\simeq&\!\!\!\!1- ( A(3Al_1/2 -H-2\rho)+3 )\times \nonumber \\
&\times&\!\!\!\!\exp\{-A(H+2\rho)\} \left[ \hat{\theta}_{l_1}(H)  \right]^{M/l_1-2}
\ee
with
\bea
\hat{\theta}_{l_1}(H)&=&\\
&&\!\!\!\!\!\!\!\!\!\!\!\!\!\!\!\!\!\!\!\!\!\!\!\!\!\!\!\!\!\!\!\!\!\!\!\!\frac{ 1- (A(3Al_1/2 -H-2\rho)+3 )\exp\{-A(H+2\rho)\}}{ 1- (A(Al_1 -H-2\rho)+3 )\exp\{-A(H+2\rho)\}} \,.
\eea
Also
\bea
\mathbb{E}_\infty \tau_{S,l_0,l_1}(H)  &\simeq& -\\
 &&\!\!\!\!\!\!\!\!\!\!\!\!\!\!\!\!\!\!\!\!\!\!\!\!\!\!\!\!\!\!\!\!\!\!\!\!\!\!\!\!\!\!\!\!\!\!\!\!\!\!\!\!\frac{l_1[ 1\!-\!(A( 3Al_1/2 \!-\!H\!-\!2\rho)\!+\!3 )\exp\{-A(H\!+\!2\rho)\}]  }{[\hat{\theta}_{l_1}(H)]^2\log(\hat{ \theta}_{l_1}(H))} \,. \nonumber
\eea
\vspace{-0.6cm}
\be\label{Approx2}
\ee

\begin{figure}[h]
\begin{center}
 \includegraphics[width=0.5\textwidth]{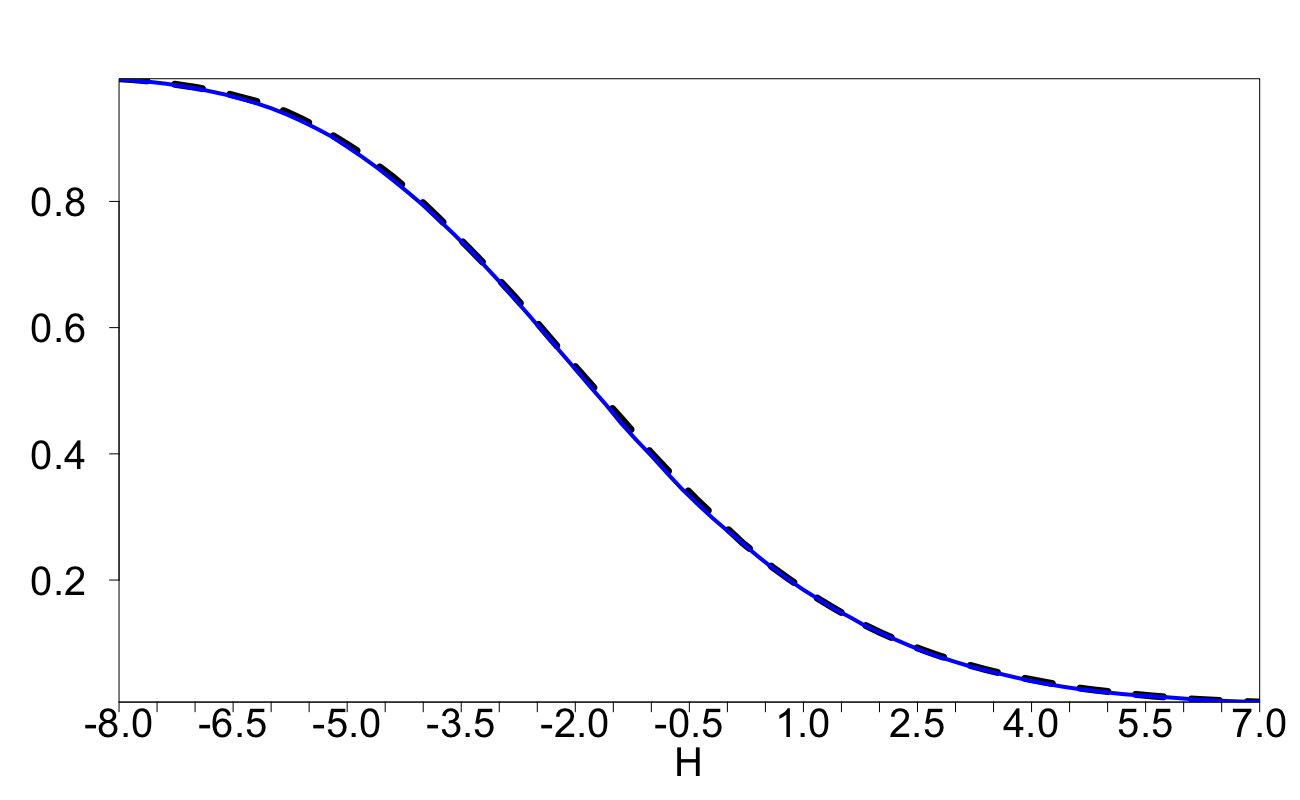}
\end{center}
\caption{Empirical probabilities of reaching the barrier $H$ (dashed black) and approximation \eqref{boundary_crossing_approx2} (solid blue):
 $A=1$, $M/l_1=4$ with $l_0=25$ and $l_1=50$.  }
\label{large_deviation0}
\end{figure}

\begin{figure}
\begin{center}
\includegraphics[width=0.5\textwidth]{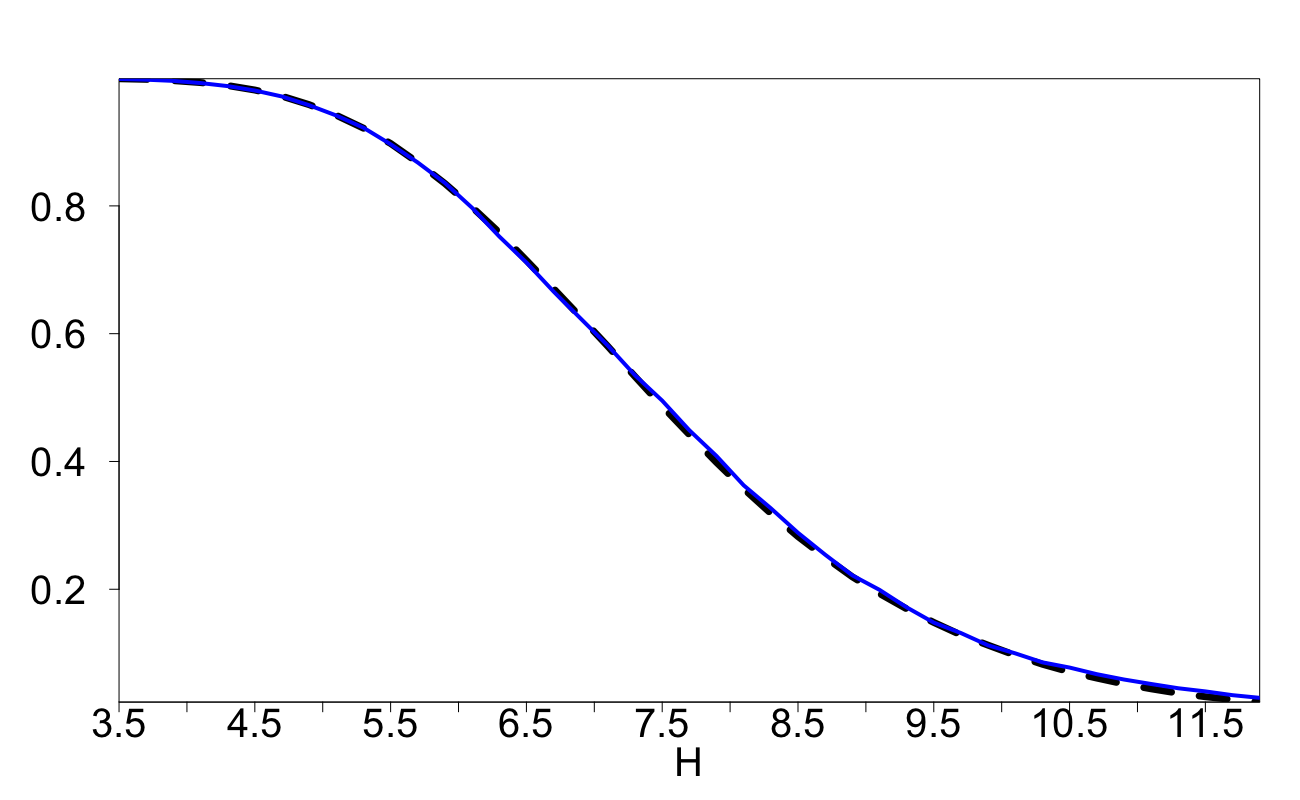}
\end{center}
\caption{Empirical probabilities of reaching the barrier $H$ (dashed black) and approximation \eqref{boundary_crossing_approx2} (solid blue): $A=0.5$, $M/l_1=25$ with $l_0=10$ and $l_1=20$.  }
\label{large_deviation00}
\end{figure}

\begin{figure}
\begin{center}
 \includegraphics[width=0.5\textwidth]{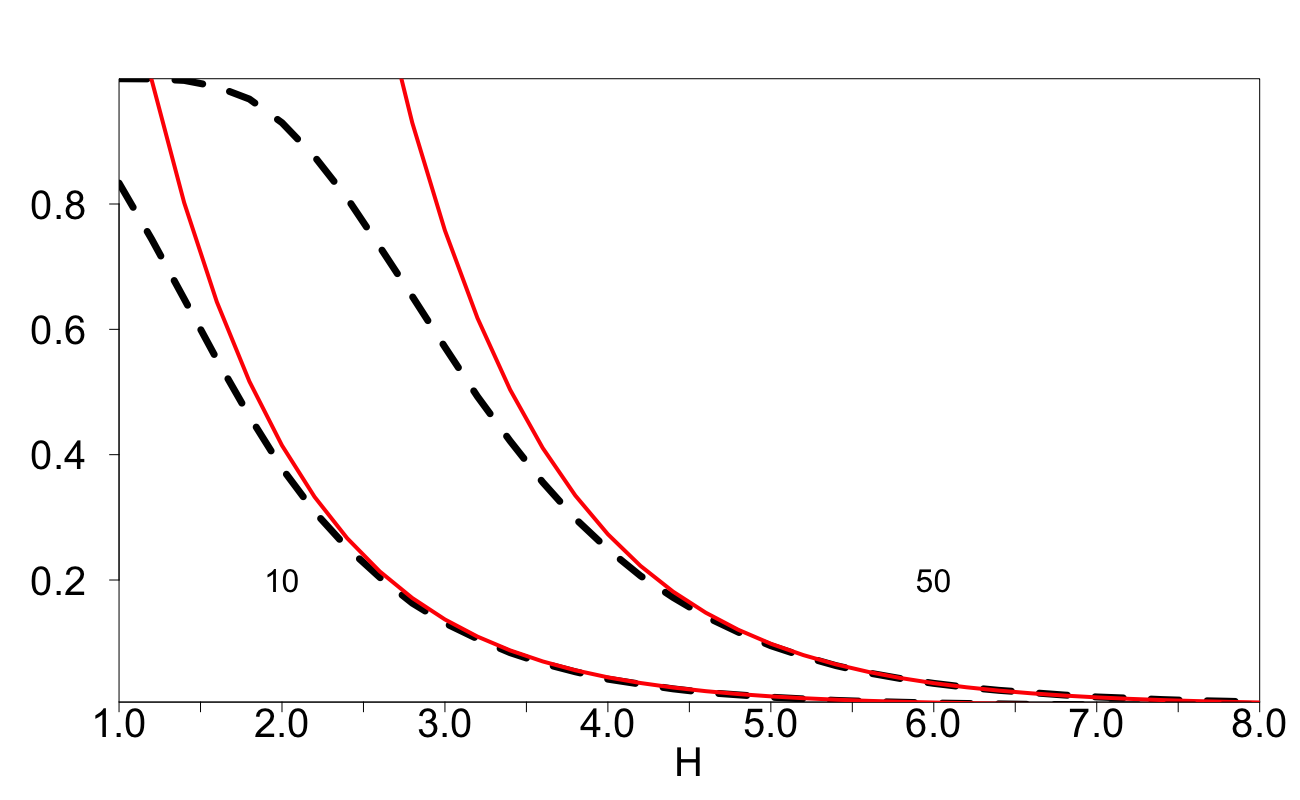}$\;\;$
\end{center}
\caption{Empirical probabilities of reaching the barrier $H$ (dashed black) and corresponding versions of approximation \eqref{Approx1} (solid red): $A=1$, $m/l_1=1$ with (a) $l_1=10$ and (b) $l_1=50$. }
\label{large_deviation11}
\end{figure}

\begin{figure}
\begin{center}
\includegraphics[width=0.5\textwidth]{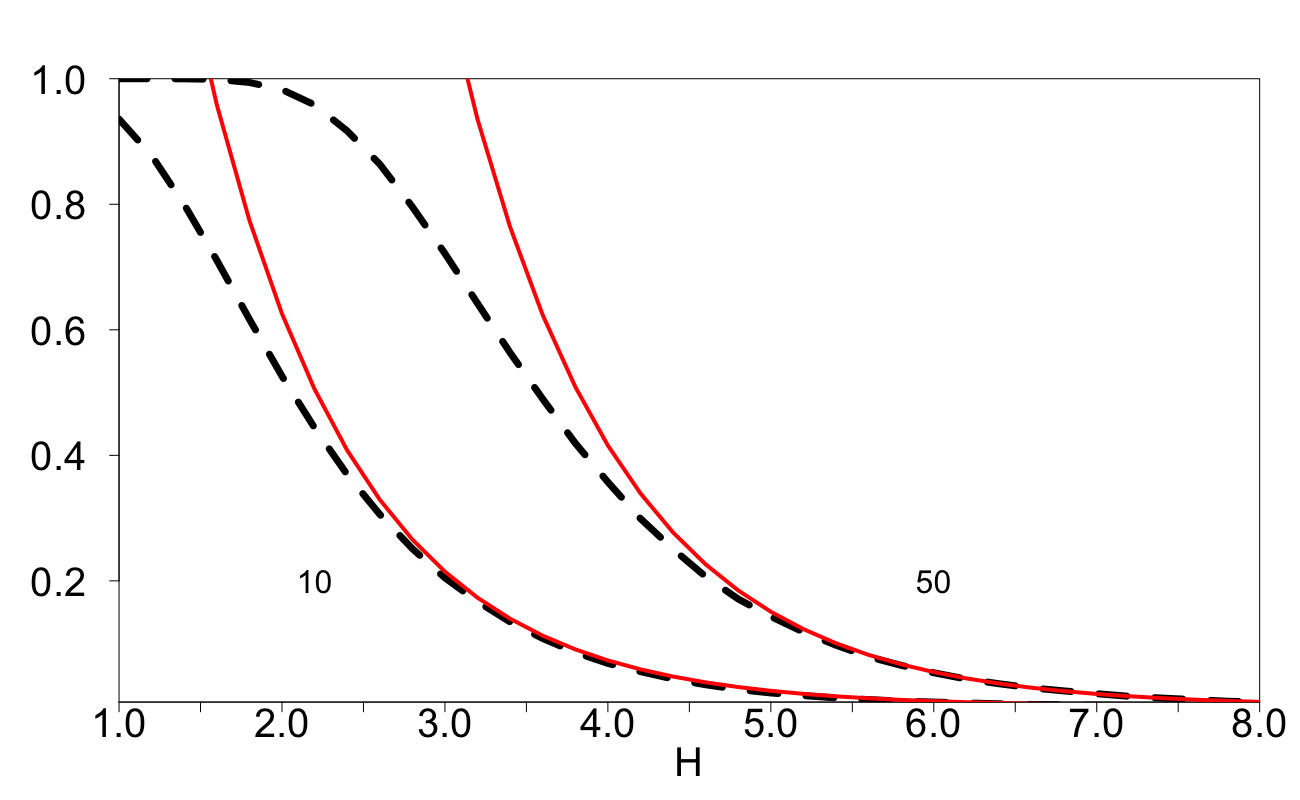}
\end{center}
\caption{Empirical probabilities of reaching the barrier $H$ (dashed black) and corresponding versions of approximation \eqref{Approx1} (solid red): $A=1$, $m/l_1=2$ with  $l_1=10$ and (b) $l_1=50$. }
\label{large_deviation1}
\end{figure}

\begin{figure}
\begin{center}
 \includegraphics[width=0.5\textwidth]{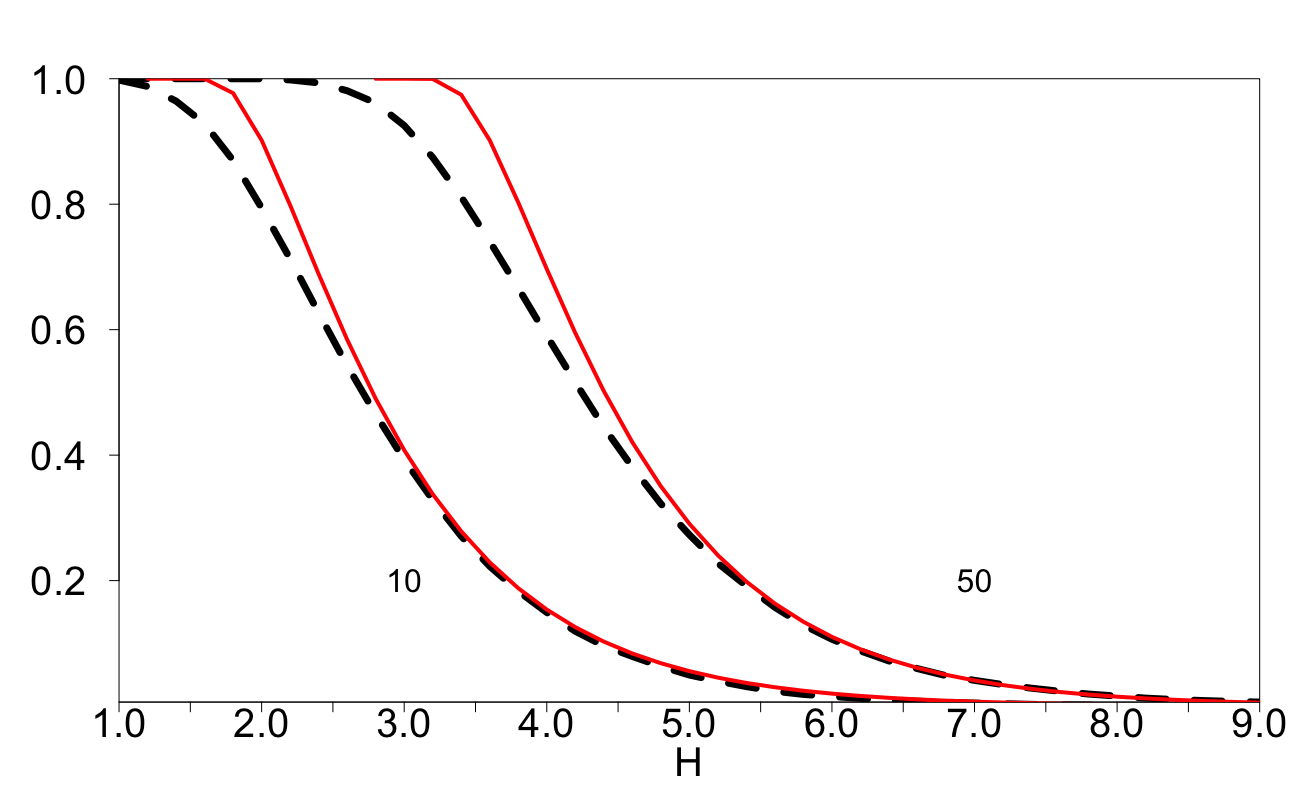}$\;\;$
\end{center}
\caption{Empirical probabilities of reaching the barrier $H$ (dashed black) and corresponding versions of approximation \eqref{Approx1} (solid red): $A=1$,  $m/l_1=10$ with  $l_1=10$ and (b) $l_1=50$. }
\label{large_deviation2}
\end{figure}

\begin{figure}
\begin{center}
\includegraphics[width=0.5\textwidth]{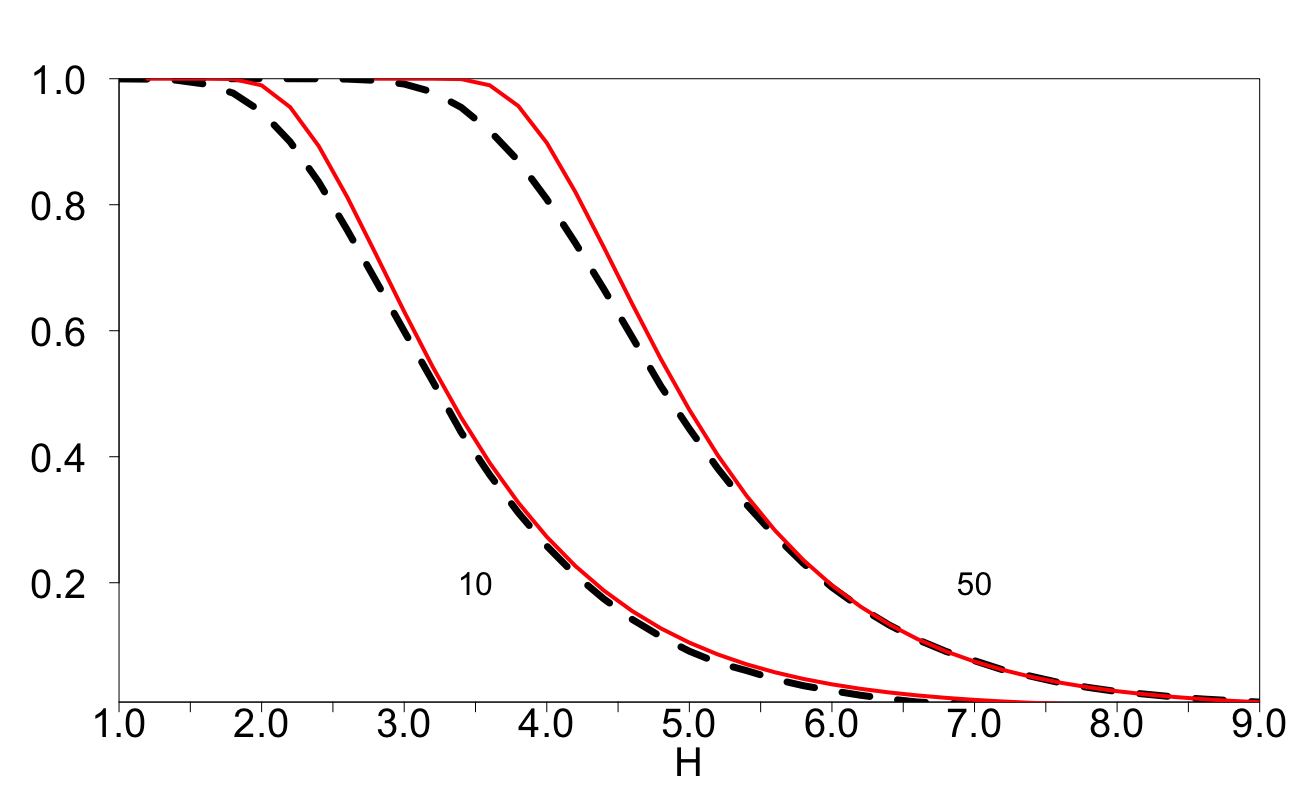}
\end{center}
\caption{Empirical probabilities of reaching the barrier $H$ (dashed black) and corresponding versions of approximation \eqref{Approx1} (solid red): $A=1$, $m/l_1=5$ with (a) $l_1=10$ and (b) $l_1=50$.}
\label{large_deviation22}
\end{figure}

The accuracy of the approximation in \eqref{boundary_crossing_approx2}  is demonstrated in Figures~\ref{large_deviation0}-\ref{large_deviation00} for different $l_0, l_1, M$ and $A$ as a function of $H$. In this approximation, ${ F}_{l_0,l_1}(H, 2l_1)$ and ${ F}_{l_0,l_1}(H, l_1)$ have been approximated using Monte Carlo simulations with 100,000 repetitions.
In these figures, the probability ${ F}_{l_0,l_1}(H; M)$ is depicted with a thick dashed black line and is obtained from simulations. The approximation in \eqref{boundary_crossing_approx2} is depicted with a solid blue line. From these figures the high accuracy of approximation \eqref{boundary_crossing_approx2} is clearly demonstrated. In Figures~\ref{large_deviation11}-\ref{large_deviation22}, we asses the accuracy of the approximation in \eqref{Approx1}. In these figures, for $A=1$ and various $M$, the probability ${ F}_{1,l_1}(H; M)$ is depicted with a thick dashed black line whereas the approximation provided in \eqref{Approx1} is shown with a solid red line. The number present on the figure is used to show the value of $l_1$ used. From these figures, we see for large $H$ the approximation in \eqref{Approx1} is adequate.
In Tables~\ref{large_deviation_table1}-\ref{large_deviation_table2} the accuracy of the approximations provided in \eqref{MOSUM_ARL2} and \eqref{Approx2} are assessed for different $H$. We see the approximation in \eqref{MOSUM_ARL2} is extremely accurate for all $H.$ For large $H$, the approximation in \eqref{Approx2} is fairly accurate and has the benefit of explicit evaluation. For small $H$ and small $A$, the accuracy of \eqref{Approx2} should deteriorate.

\begin{table*}
\caption{Approximations for $\mathbb{E}_\infty \tau_{S,l_0,l_1}(H)$  with $l_0=25, l_1=50, A=1$.}
\label{large_deviation_table1}
\centering
\begin{tabular*}{\textwidth}{@{\extracolsep{\fill}} c c c c c c c c }
\hline
$H$ & -5 & -4.5 & -4&-3.5 & -3 &-2.5 & -2  \\
\hline
  \eqref{MOSUM_ARL2}& 78 & 96   & 116  & 144 &  177  & 218 &  271          \\
$\mathbb{E}_\infty \tau_{S,l_0,l_1}(H)$  & 77 & 95 & 115  &144  & 179 & 217  &276      \\
\hline
\end{tabular*}
\end{table*}

\begin{table*}
\caption{Approximations for $\mathbb{E}_\infty \tau_{S,1,l_1}(H)$  with $l_1=10, A=1$.}
\label{large_deviation_table2}
\centering
\begin{tabular*}{\textwidth}{@{\extracolsep{\fill}} c c c c c c c c }
\hline
$H$ & 2 & 2.25 &2.5 &2.75 & 3 &3.25 & 3.5  \\
\hline
 \eqref{Approx2}& 20 &32   &49  & 71  & 100   & 137 &  185          \\
  \eqref{MOSUM_ARL2}& 29 &43   &58  & 82  & 109   & 147 &  201          \\
$\mathbb{E}_\infty \tau_{S,l_0,l_1}(H)$  & 30 & 43 & 59  &81  & 111 & 148  &201      \\
\hline
\end{tabular*}
\end{table*}

\subsubsection{The presence of nuisance parameters}
  
Here we briefly consider statistics aimed at detecting a transient change when certain nuisance parameters require estimation. The brevity of this discussion is because in practice for online change-point problems, the behaviour of the time series under the null hypothesis of no change-point is often observed for a lengthy period of time. This allows for the accurate estimation of certain nuisance parameters and they can therefore be assumed known. Many of the following statistics appear in some form in \cite{yao1993tests} when addressing the offline change-point problem, and a number of approximations for the false alarm error are provided.
The log likelihood ratio given in \eqref{MOSUM_likelihood}, where $f$ and $g$ are given in \eqref{Gaussian_problem2}, is
 \bea
\Lambda_{\nu,\nu+l}=A \sum_{j=\nu+1}^{n+l} \left( y_j-\mu-\frac{A}{2} \right) \,.
 \eea
Using motivation from \cite{levin1985cusum}, if $\mu$ is unknown, $l$ is unknown but bounded $l_0\leq l\leq l_1$ and $A$ is known, then one can replace $\mu$ with its maximum likelihood estimator under $H_\infty$; $\hat{\mu} := \sum_{i=1}^ny_i/n$ to obtain:
\bea
Z^1_{n}:=\max_{\substack{0\leq \nu <\nu+l \leq n\\l_0 \leq l \leq l_1} } A \sum_{j=\nu+1}^{\nu+l} \left( y_j-\hat{\mu}-\frac{A}{2} \right)\,.
\eea

In \cite{siegmund1986boundary}, $\mu$ was replaced with its average over the null and alternative hypotheses to obtain the true likelihood ratio statistic:
\bea
Z^2_{n}:=\max_{\substack{0\leq \nu <\nu+l \leq n\\l_0 \leq l \leq l_1} } A \sum_{j=\nu+1}^{\nu+l} \left( y_j-\hat{\mu}-\frac{A}{2}\left(1-\frac{l}{n} \right) \right).
\eea

If $\mu$ and $A$ are both unknown, the square root of the log likelihood ratio statistic is:
\bea
Z^3_{n}:=\max_{\substack{0\leq \nu <\nu+l \leq n\\l_0 \leq l \leq l_1} }  \frac{\left(\sum_{j=\nu+1}^{\nu+l} y_j\right) -l \hat{\mu} }{\sqrt{l(1-\frac{l}{n})}} \,.
\eea

It is not obvious how one can translate the offline change-point results of \cite{levin1985cusum,siegmund1986boundary,hogan1986large,siegmund1988approximate,yao1989large} to address the online change point problem in the presence of nuisance parameters. This is because the change-point statistics $Z^1_{n}, Z^2_{n}$ and $Z^3_{n}$ can no longer be written recursively as  the estimators for the unknown parameters get updated at each time $n$.

%\begin{table}[h]
%\centering
%\begin{tabular}{|c||c|c|c|c|c|c|c|}
%\hline
%$H$ & 2 & 2.25 &2.5 &2.75 & 3 &3.25 & 3.5   \\
%\hline
% \eqref{Approx2}&  471 & 791  & 1392   & 2587  & 5099    &  10695  &    23918           \\
%$\mathbb{E}_\infty \tau_S(H) $ & 472  &792  & 1397  &  2588 & 5085  & 10749   &    24131      \\
%\hline
%\end{tabular}
%\caption{Approximations for $\mathbb{E}_\infty \tau_S(H)$ with $L=50$.}
%\label{large_deviation_table2}
%\end{table}
%

\subsection{Optimality criteria} \label{Comparing_procedures2}

 For online detection of transient changes, optimality criteria like \eqref{Pollak_sad} and \eqref{lorden} do not have much meaning as the change in distributions is not permanent (signal can be missed). Instead, optimality involving the maximisation of the probability of detection under a constraint on the false alarm risk is more applicable, see \cite{han1999some,bakhache2000reliable}. One could use a worst-case criterion of the form:
\be\label{Power1}
&& \inf_{\nu  }{\rm P}_{\nu}\{\tau(H) -\nu <\! T  \, | \, \tau(H) \!>\! \nu  \} \,,
\ee
where $T>1$ is the maximum length of time after the change-point occurs that it must be detected; this is problem specific and is therefore chosen by the user.  By imposing the condition of a long run with no false alarms, another possible criterion is
\be\label{Power2}
&&  \lim_{\nu \rightarrow \infty}{\rm P}_{\nu}\{\tau(H) -\nu <\! T  \, | \, \tau(H) \!>\! \nu  \}.
\ee

Using ARL as the measure of false alarm risk, a stopping rule $\tau \in \Delta(C)$ is then optimal for a given $C$ if it maximises \eqref{Power1} or \eqref{Power2}; recall $\Delta(C) = \{\tau: E_{\infty}\tau \geq C \},\,\,\, C>1.$

\subsubsection{MOSUM procedure}
For the MOSUM procedure given in \eqref{MOSUM_normal}, the quantity \eqref{Power2} was the focus of study in \cite{noonan2020power} and built on the continuous time results of \cite{zhigljavsky2021first}.
For $T=l+L$, the quantity  \eqref{Power2} is equivalent to: 
\be\label{Power3}
 {\cal{P}}_S(H,A,L)\! &:=&\!\!\! \lim_{\nu \rightarrow \infty}{\rm P}_{\nu}\{S_{n,L}\!>\!H \nonumber \\
&& \!\!\!\!\!\!\!\!\!\!\!\!\!\!\!\!\!\!\!\!\!\!\!\!\!\!\!\!\!\!\!\!\!\!\!\! \text{  for some }  n\!\in\![\nu'\!+\!1,\nu\!+\!l\!-\!1  ]  \, | \, \tau_{S,L}(H) \!>\! \nu'  \}, \nonumber \\
\ee
with $\nu':=\nu-L$.

Formally, we require   $\nu \to \infty$ in \eqref{Power3}. This is to ensure that the sequence of moving sums $\{S_{n,L}\}_n$ reaches the stationary behaviour under the null hypothesis and given that we have not crossed the threshold $H$.  However, as discussed \cite{noonan2020power,zhigljavsky2021first}, this stationary regime is reached very quickly and in all approximations  it is  enough to only require  $\nu \geq 2L$.

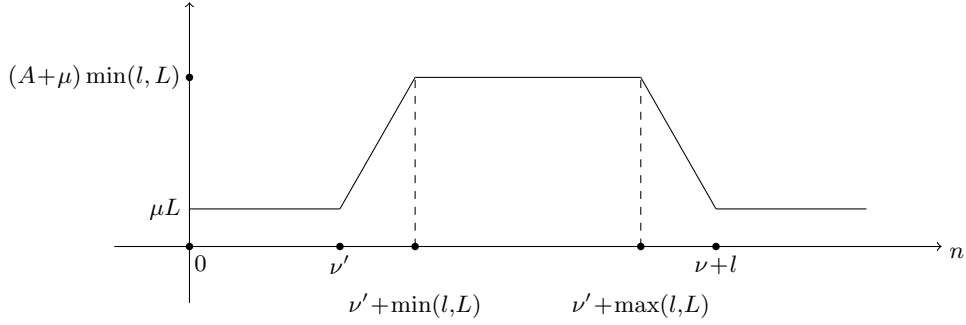
\begin{figure*}
\begin{center}
\begin{tikzpicture}

      \draw[->] (-4,-0.25) -- (7,-0.25) node[right] {};
      \draw[->] (-3,-1) -- (-3,3) node[right] {};
      \draw[scale=1,domain=-3:-1,smooth,variable=\x] plot ({\x},{0.25});
      \draw[scale=1,domain=-1:0,smooth,variable=\x] plot ({\x},{0.25+1.75*(\x+1)});
       \draw[scale=1,domain=0:3,smooth,variable=\x] plot ({\x},{2});
      \draw[scale=1,domain=4:6,smooth,variable=\x] plot ({\x},{0.25});
      \draw[scale=1,domain=3:4,smooth,variable=\x] plot ({\x},{2-1.75*(\x-3)});
      \node at (-3,-0.25) [circle,fill,inner sep=1pt]{};
      \draw [-] (-2.85,-0.25) node [below] {0};
      \draw [-] (-3,0.25) node [left] {$\mu L$};
      \draw [-] (-3,2) node [left] {$(A\!+\!\mu)\min(l,L)$};
      \node at (-3,2) [circle,fill,inner sep=1pt]{};
      \node at (-1,-0.25) [circle,fill,inner sep=1pt]{};
      \draw [-] (-1,-0.25) node [below] {$\nu'$};
      \node at (4,-0.25) [circle,fill,inner sep=1pt]{};
      \draw [-] (4,-0.25) node [below] {$\nu \!+\!l$};

      \draw [-] (3,-0.75) node [below] {$\nu'\!+\!\max(l,\!L)$};
      \node at (3,-0.25) [circle,fill,inner sep=1pt]{};

      \draw [-] (0,-0.75) node [below] {$\nu'\!+\!\min(l,\!L)$};
      \node at (0,-0.25) [circle,fill,inner sep=1pt]{};
      \draw[black,dashed] (0,-0.25) -- (0,2);
      \draw[black,dashed] (3,-0.25) -- (3,2);
       \draw [-] (7.2,-0.15) node [below] {$n$};

  \end{tikzpicture}
 \end{center}
\vspace{-0.3cm} \caption{$\mathbb{E}_\nu S_{n,L}$ as a function of $n$}
 \label{means_of_S_nA}
 \end{figure*}

The reasoning behind the choice $T=l+L$ is as follows.
Assume $\mathbb{H}_\nu$ with $\nu < \infty$, and that $\nu$ is suitably large. If the barrier $H$ is reached for any sum $S_{n,L}$ with $n\leq \nu'$ then, since  there are no parts of the signal in the sums $S_{0,L}, \ldots S_{\nu',L}$, we classify  the event of reaching the barrier as a  false alarm.
Each one of the sums $S_{\nu'+1,L}, \ldots, S_{\nu+l-1,L}$ has mean larger than $L\mu$ as it  contains at least a part of the signal.
Reaching the barrier $H$ by any of these sums will be classified as a correct detection of the signal. If neither of these sums reaches $H$, then
 we say that we failed to detect the signal and further events when $S_{n,L} \geq H $ with $n\ge \nu+l $
will again be classified as  false alarms.
In Figure~\ref{means_of_S_nA} we display the values $\mathbb{E}_\nu S_{n,L}$ as a function of $n$.

Define the function
\bea
Q(n;A,L,\nu')\!:=\!\left\{ \begin{array}{cl}
 0 &  \\
     A(n-\nu') &  \\
      A\min(l,L)&   \\
     \!\!\!\!  A(\min(l,L)\!-\!(n\!-\!\nu'\!-\!\max(l,L))&  
       \end{array} \right.
\eea
\bea
\begin{array}{cl}
  & \text{   for } n \le \nu'\, \text{ or }  n \ge \nu+l\,  \\
     & \text{   for } \nu'<n\le \nu'+\min(l,L)\,  \\
& \text{   for }   \nu'\!+\!\min(l,L)<n\le \nu'\!+\!\max(l,L) \, \\
     & \text{   for }   \nu'+\max(l,L)<n\le \nu+l-1 \,.
       \end{array} 
\eea
By subtracting $\mathbb{E}_\nu S_{n,L}$ from the threshold $H$ and standardising the random variables $S_{n,L}$ the power of the test given in \eqref{Power3} can be expressed in terms of probability under $\mathbb{H}_\infty$:
\be\label{power4}
\!{\cal P}_\xi(h,A,L) \!:=\!\!  \lim_{\nu \rightarrow \infty}\!{\rm P}_\infty \! \bigg\{\xi_{n,L}\!>\!h\! -\! \frac{Q(n;A,L,\nu')}{\sigma\sqrt{L}} \nonumber \\
\text{ for some } n\!\in\![\nu'\!+\!1,\nu\!+\!l\!-\!1  ] \,  \big| \, {\tau}_{\xi}(h) \!>\! \nu'  \bigg \}  \!, \nonumber \\
\ee
where ${\cal{P}}_S(H,A,L)={\cal P}_\xi(h,A,L)$.
To approximate ${\cal P}_\xi(h,A,L)$, the approach taken in \cite{noonan2020power} was similar to the approach take to approximate ARL. The approach is as follows. We firstly approximate the problem in the continuous-time setting and compute probabilities for the Gaussian process $S(t)$. Then, use the results of D. Siegmund to correct the continuous time probability for discrete time. Fix $\gamma=A\sqrt{L}/\sigma$, $\kappa=\nu' /L$, $\lambda=l/L$ and define the function
\bea
Q(t;\gamma,\kappa,\lambda)=\left\{ \begin{array}{cl}
 0 &   \\
     \gamma(t-\kappa) &  \\
      \gamma\min(1,\lambda)&  \\
       \gamma(\min(1,\lambda)-(t-\kappa-\max(1,\lambda))&  \,
       \end{array} \right.
\eea
\bea
\begin{array}{cl}
    & \text{   for } t \le \kappa \text{ or } t \ge \kappa+1+\lambda. \, \\
     & \text{   for } \kappa<t\le\kappa+\min(1,\lambda)\,  \\
     & \text{   for }   \kappa+\min(1,\gamma)<t\le \kappa+\max(1,\lambda) \, \\
   & \text{   for }   \kappa+\max(1,\lambda)<t\le \kappa+1+\lambda \,.
       \end{array} 
\eea

The diffusion approximation for the power of the test is
\be\label{Diffusion_power}
{\cal P}(h,A):=  \lim_{\kappa \rightarrow \infty}{\rm P}_\infty\{S(t)>h-Q(t;\gamma,\kappa,\lambda) \nonumber \\ \text{  for some  }  n\in[\kappa,\kappa+1+\lambda  ]  \, | \, \tilde{\tau}(h) > \kappa \},\nonumber \\
\ee
where $\tilde{\tau}(h)=\inf\{t>0: S(t)>h\}$. That is, we make the approximation
\bea \label{diffusion_approx}
{\cal P}_\xi(h,A,L) \cong {\cal P}(h,A)
\eea
by assuming $L \to \infty$.

The complexity of  computation of  the diffusion approximation $ {\cal P}(h,A)$ and its discrete-time corrected version  depends on the choice of $L$ in comparison to $l$.
Here, we will only consider the scenario of $\lambda=l/L=1$ which corresponds to the case of $l$ known at the MOSUM construction stage. The two other cases of 
$\lambda >1$ and $\lambda <1$ are studied in \cite{noonan2020power}.

For $\lambda=1$, the diffusion approximation for ${\cal P}_\xi(h,A,L)$  given in \eqref{Diffusion_power} reduces to
\be\label{Case_1_diffusion}
{\cal P}(h,A) \!
=\! \lim_{\kappa \rightarrow \infty} {\rm P}_\infty \{ S(t) \geq
h -Q(t;\gamma,\kappa)\,\, \nonumber\\ \mbox{for some $t \!\in$ \!$[\kappa, \kappa +2]$}  \big | \,   \tilde{\tau}(h) > \kappa    \}\! , \nonumber \\
\ee
where
 $Q(t;\gamma,\kappa) = \gamma \max \left \{  0, 1-|t-(\kappa+1) |  \right \}$. The barrier $h-Q(t;\gamma,\kappa)$ is depicted in Figure~\ref{two_cases}.

\begin{figure*}
\begin{minipage}[b]{0.5\linewidth}
\centering
\begin{tikzpicture}[scale=0.8]

      \draw[->] (-1.5,0) -- (5,0) node[right] {};
      \draw[->] (0,-2.5) -- (0,4.2) node[above] {};
      \draw[scale=1,domain=-1.5:0,smooth,variable=\x] plot ({\x},{3});
      \draw[scale=1,domain=0:2,smooth,variable=\x]  plot ({\x},{3-2.5*\x});
      \draw[scale=1,domain=2:4,smooth,variable=\x]  plot ({\x},{-2+2.5*(\x-2)});
      \draw[scale=1,domain=4:4,smooth,variable=\x] plot ({\x},{3});
       \draw [-] (-0.2,3) node [above] {$h$};
       \draw [-] (2,0) node [above] {$\kappa\!+\!1$};
       \draw [-] (4,0) node [above] {$\kappa\!+\!2$};
       \draw [-] (0.2,0) node [above] {$\kappa$};
       \draw [-] (0.6,1.5) node [right, font=\small] {$h\!-\!\gamma(t\!-\!\kappa) $};
       \draw [-] (3.4,1.5) node [right,font=\small] {$(h\!-\!\gamma)+$};
       \draw [-] (3.4,1.1) node [right,font=\small] {$\!\gamma(t\!-\!\kappa\!-\!1)$};
       \draw [-] (0,-2) node [left] {$h\!-\!\gamma$};
       \node at (2,0) [circle,fill,inner sep=1pt]{};
       \node at (4,0) [circle,fill,inner sep=1pt]{};
       \node at (0,0) [circle,fill,inner sep=1pt]{};
       \node at (0,-2) [circle,fill,inner sep=1pt]{};
        \draw [-] (5,-0.15) node [below] {$t$};

  \end{tikzpicture}
\caption{ Barrier $h-Q(t;\gamma,\kappa)$ for $\lambda=1$.}
\label{two_cases}
\end{minipage}
\hspace{-0.2cm}
\begin{minipage}[b]{0.45\linewidth}
\centering
\begin{tikzpicture}[scale=0.8] %[x={10.0pt},y={10.0pt}]
      \draw[->] (0,0) -- (7,0) node[right] {};
      \draw[->] (0,-2.5) -- (0,4.2) node[above] {};
      \draw[scale=1,domain=0:2,smooth,variable=\x]  plot ({\x},{3});
      \draw[scale=1,domain=2:4,smooth,variable=\x]  plot ({\x},{3-2.5*(\x-2)});
       \draw[scale=1,domain=4:6,smooth,variable=\x]  plot ({\x},{-2+2.5*(\x-4)});
       \draw [-] (-0.2,3) node [above] {$h$};
       \draw [-] (2,0) node [above] {$1$};
       \draw [-] (4,0) node [above] {$ 2$};
        \draw [-] (6,0) node [above] {$ 3$};
       \draw [-] (0.2,0) node [above] {$0 $};
       \draw [-] (2.7,1.5) node [right, font=\small] {$h\!+\!\gamma\!-\!\gamma t $};
       \draw [-] (5.4,1.5) node [right, font=\small] {$h\!-\!3\gamma\!+\!\gamma t $};
       \draw [-] (0,-2) node [left] {$h\!-\!\gamma$};
       \node at (2,0) [circle,fill,inner sep=1pt]{};
         \node at (6,0) [circle,fill,inner sep=1pt]{};
       \node at (4,0) [circle,fill,inner sep=1pt]{};
       \node at (0,0) [circle,fill,inner sep=1pt]{};
       \node at (0,-2) [circle,fill,inner sep=1pt]{};
        \draw [-] (7,-0.15) node [below] {$t$};
\end{tikzpicture}
\caption{Barrier $B(t;h,0,-\gamma,\gamma)$.}
\label{three_cases}
\end{minipage}\end{figure*}

The probability \eqref{Case_1_diffusion} was considered in \cite{zhigljavsky2021first}, where approximations accurate to more than $4$ decimal places were developed.
Define the following two conditional probabilities:
\bea\label{piecewise-problem3}
 F_{h,0}(1|x) \!\!\!\!&:=&\!\!\!\! {\rm P}_\infty(S(t)< h \text{  for all   } t\in[0,1]\,\, \\
 && | \,\,S(0)=x) \,,\\
F_{h,0,-\gamma,\gamma}(3 | x) \!\!\!\!&:=&\!\!\!\! {\rm P}_\infty(S(t)< B(t;h,0,-\gamma,\gamma) \\
&& \text{  for all   } t\in[0,3]\,\, | \,\,S(0)=x),
\eea
where the barrier $B(t;h,0,-\gamma,\gamma)$ is defined as
\bea
B(t;h,0,-\gamma,\gamma)=
\left \{ \begin{array}{ll}
\!\!\! h, &   0 \le t \le 1\\
\!\!\!h-\gamma (t-1), &   1 < t \le 2\\
\!\!\!h-\gamma +\gamma (t-2), &  2 <t\le 3 \\
\!\!\!0 &   \mbox{otherwise,}
\end{array} \right.
\eea
 and is depicted in Figure~\ref{three_cases}. From \cite{zhigljavsky2021first} we obtain
\be\label{Approximation1}
{\cal P}(h,A) \cong1- \frac{F_{h,0,-\gamma,\gamma}(3 | 0)
} { F_{h,0}(1|0) }\,,
\ee
where
\be \label{F_h_1}
&&  F_{h,0}(1\, |\, x)=\Phi(h)-\exp\left(-(h^2-x^2)/2 \right)\Phi(x) \nonumber  \\
 \ee
 and
 \bea
&&\!\!\!\!\!\!\!\!\!F_{h,0,-\gamma,\gamma}(3\, |\, x)=\frac{e^{{\gamma}^2/2}}{\varphi(x)} \int_{-x-h}^{\infty}\int_{x_2-h+\gamma}^{\infty} e^{ -\gamma(x_3-x_2)}    \nonumber \\
&&\!\!\!\! \times \det \!\! \left[ \begin{matrix}
    \varphi(x)      &  \varphi(-x_2\!-\!h)  \\
    \varphi(h)      & \varphi(-x\!-\!x_2) \\
    \varphi(x_2\!+\!2h\!+\!x)      & \varphi(h)  \\
        \varphi(x_3\!+\!3h\!-\!\gamma\!+\!x)      \!\!& \varphi(x_3\!+\!2h\!-\!\gamma\!-\!x_2) \\
\end{matrix}\!\! \right. \nonumber\\
&&\ \!\! \left. \begin{matrix}
   \varphi(-x_3\!-\!2h\!+\!\gamma) &  \Phi(-x_3\!-\!2h\!+\!\gamma) \\
  \varphi(-x\!-\!x_3\!-\!h\!+\!\gamma)\!\! & \!\Phi(-x\!-\!x_3\!-\!h\!+\!\gamma) \\
   \varphi(x_2\!-\!x_3\!+\!\gamma)& \Phi(x_2\!-\!x_3\!+\!\gamma) \\
   \varphi(h)& \Phi(h) \\
\end{matrix}\right] \!\!dx_3dx_2 \,.\nonumber
\eea

To compute the approximation \eqref{Approximation1} one needs to numerically evaluate a two-dimensional integral which is a routine problem for modern computers.

Correcting approximation \eqref{Approximation1} for discrete time can be performed in the same manner as correcting the ARL approximations in Section~\ref{sec:Mosum}. This results in the approximation
\be
\label{eq:appr1}
{\cal P}_\xi(h,A,L)  \cong 1- \frac{F_{h_L,0,-\gamma,\gamma}(3 | 0)
} { F_{h_L,0}(1|0) }\,, \nonumber \\ \text{ where } h_L := h +\omega_L. 
\ee

\begin{figure}
\begin{center}
 \includegraphics[width=0.5\textwidth]{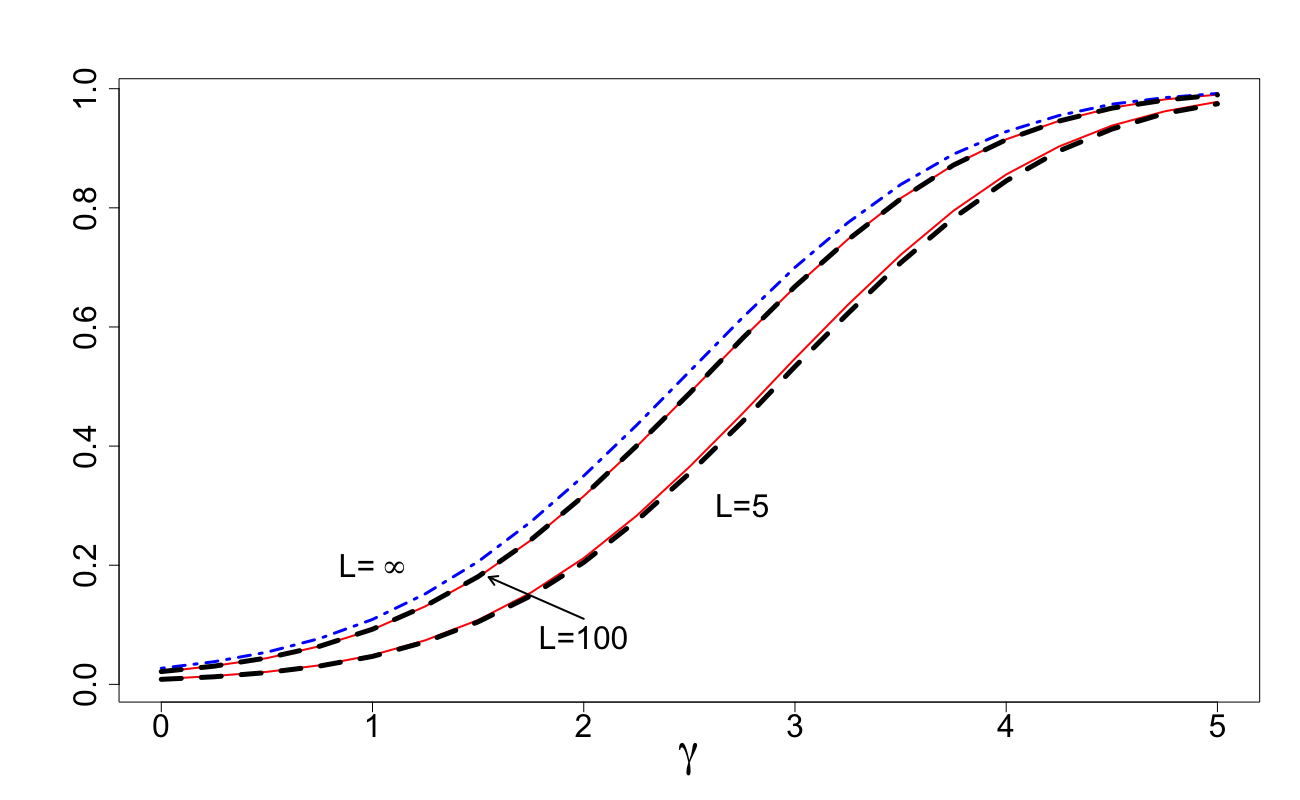}
\end{center}
\caption{Empirical probabilities of ${\cal P}_\xi(h,A,L)$ (thick dashed black) and its approximations (solid red and solid blue) for $h=3$.}
\label{fig:case_one}
\end{figure}

\begin{figure}
\begin{center}
\includegraphics[width=0.5\textwidth]{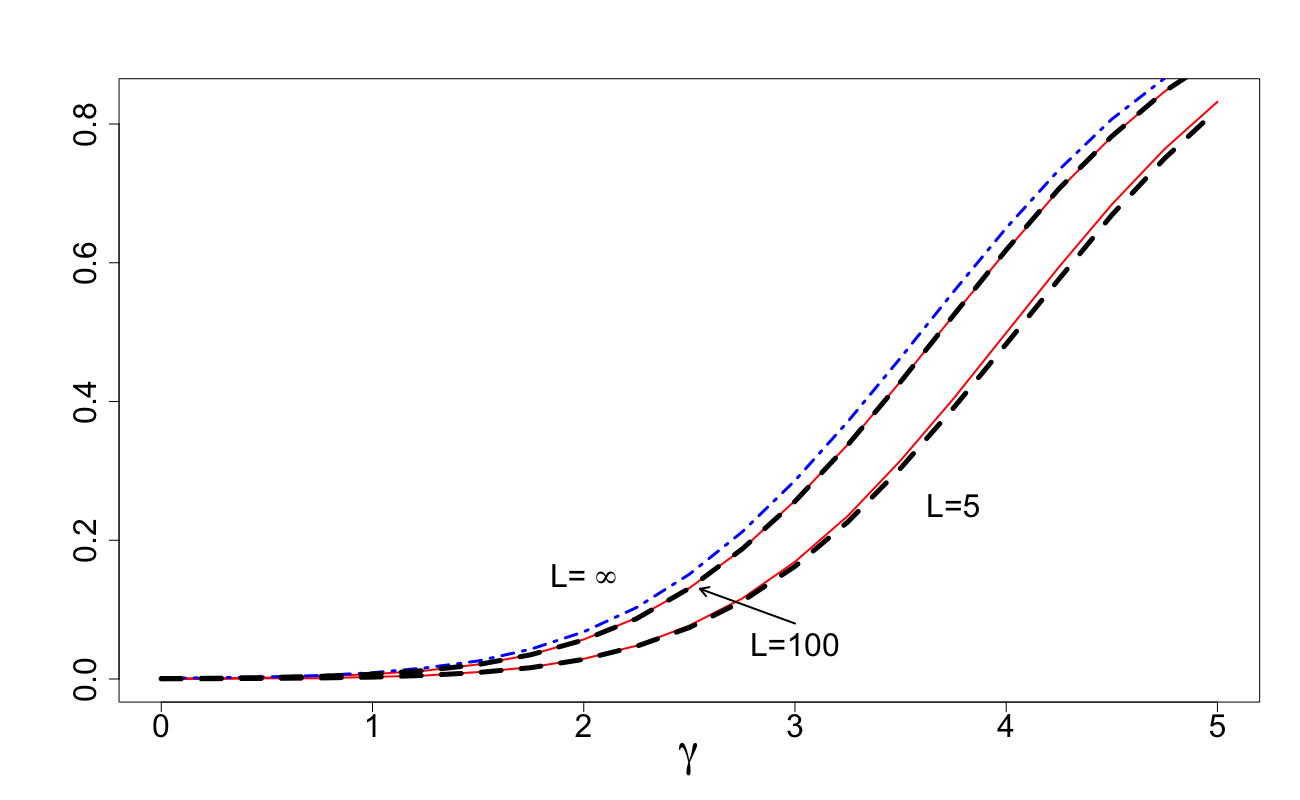}
\end{center}
\caption{Empirical probabilities of ${\cal P}_\xi(h,A,L)$ (thick dashed black) and its approximations (solid red and solid blue) for $h=4$.}
\label{fig:case_one1}
\end{figure}

 In Figures~\ref{fig:case_one}-\ref{fig:case_one1}, the thicker black dashed line corresponds to the empirical values of the BCP ${\cal P}_\xi(h,A,L) $ computed from 100\,000 simulations  with different values of $L$ and $\gamma$, where $\mu=0$ and $\sigma=1$. The solid red line corresponds to the approximation in \eqref{eq:appr1}. The dot-dashed blue line corresponds to the diffusion approximation given in \eqref{Approximation1}. The axis are: the $x$-axis shows the value of $\gamma$. The $y$-axis denotes the probabilities of reaching the barrier. The graphs, therefore,  show the empirical probabilities of ${\cal P}_\xi(h,A,L) $ and values of approximation \eqref{eq:appr1}.

From Figures~\ref{fig:case_one}-\ref{fig:case_one1}, we see that approximation \eqref{eq:appr1} is very accurate even for a very small  $L=5$. We also see the significance of the discrete-time correction; whilst the diffusion approximation provides sensible results should you compare it with $L=100$, for $L=5$ the diffusion approximation is very far off.

\subsection{Comparison of tests}
In this section, we compare the power of the MOSUM test in \eqref{MOSUM_normal} against the generalised MOSUM statistic \eqref{generalised_MOSUM} and the CUSUM test given in \eqref{CUSUM_test} specialised for this Gaussian example when used to detect a transient change.
% a detailed description of the CUSUM test including properties such as ARL can be found in \cite[Section 8.2]{tartakovsky2014sequential}. 
Secondly, but also simultaneously, we compare the power of the MOSUM test as $\lambda=l/L$ varies in $[0.5,2]$; the purpose is to demonstrate when the generalised MOSUM statistic becomes beneficial when the exact value of $l$ is unknown. Here, we shall consider the power criterion given in \eqref{Power2} and set $T=2l$. For the MOSUM test, the power is then
\bea
{\cal{P}}_S(H_1,A,L)\! :=\!\!\! \lim_{\nu \rightarrow \infty}{\rm P}_{\nu}\{S_{n,L}\!>\!H_1  \text{  for some } \\  n\!\in\![\nu-L\!+\!1,\nu\!-L+\!2l\!-\!1  ]  \, | \, \tau_{S,L}(H_1) \!>\! \nu-L  \}.
\eea
For the generalised MOSUM test, the power is
\bea
{\cal{P}}_Z(H_2,A,l_0,l_1)\! :=\!\!\! \lim_{\nu \rightarrow \infty}{\rm P}_{\nu}\{{S}_{n,{l_0,l_1}}\!>\!H_2  \text{  for some } \\ n\!\in\![\nu-l_1\!+\!1,\nu\!-l_1+\!2l\!-\!1  ]  \, | \, \tau_{S,l_0,l_1}(H_2) \!>\! \nu-l_1  \}.
\eea
The power of the CUSUM test for the transient change considered in then equivalent to
\bea
{\cal{P}}_{V}(H_{3},A)\! :=\!\!\! \lim_{\nu \rightarrow \infty}{\rm Pr}_{\nu}\{V_{n}\!>\!H_{3}  \text{  for some } \\ n\!\in\![\nu+\!1,\nu\!+\!2l\!-\!1  ]  \, | \, \tau_{V}(H_{3}) \!>\! \nu  \}.
\eea
To compare the three tests, the thresholds $H_1$, $H_2$ and $H_3$ have been set such that $\mathbb{E}_\infty\tau_M(H_1)=\mathbb{E}_\infty\tau_{Z}(H_{2})=\mathbb{E}_\infty\tau_{V}(H_{3})=500$. Determination of $H_{1}$ for MOSUM has been computed using the accurate approximation in \eqref{MOSUM_ARL}. For the generalised MOSUM procedure, $H_2$ is found via Monte Carlo simulations. Determination of $H_{3}$ for CUSUM was obtained using tabulated values given in \cite[p. 3237]{moustakides2009numerical}.

\begin{figure}
\begin{center}
 \includegraphics[width=0.5\textwidth]{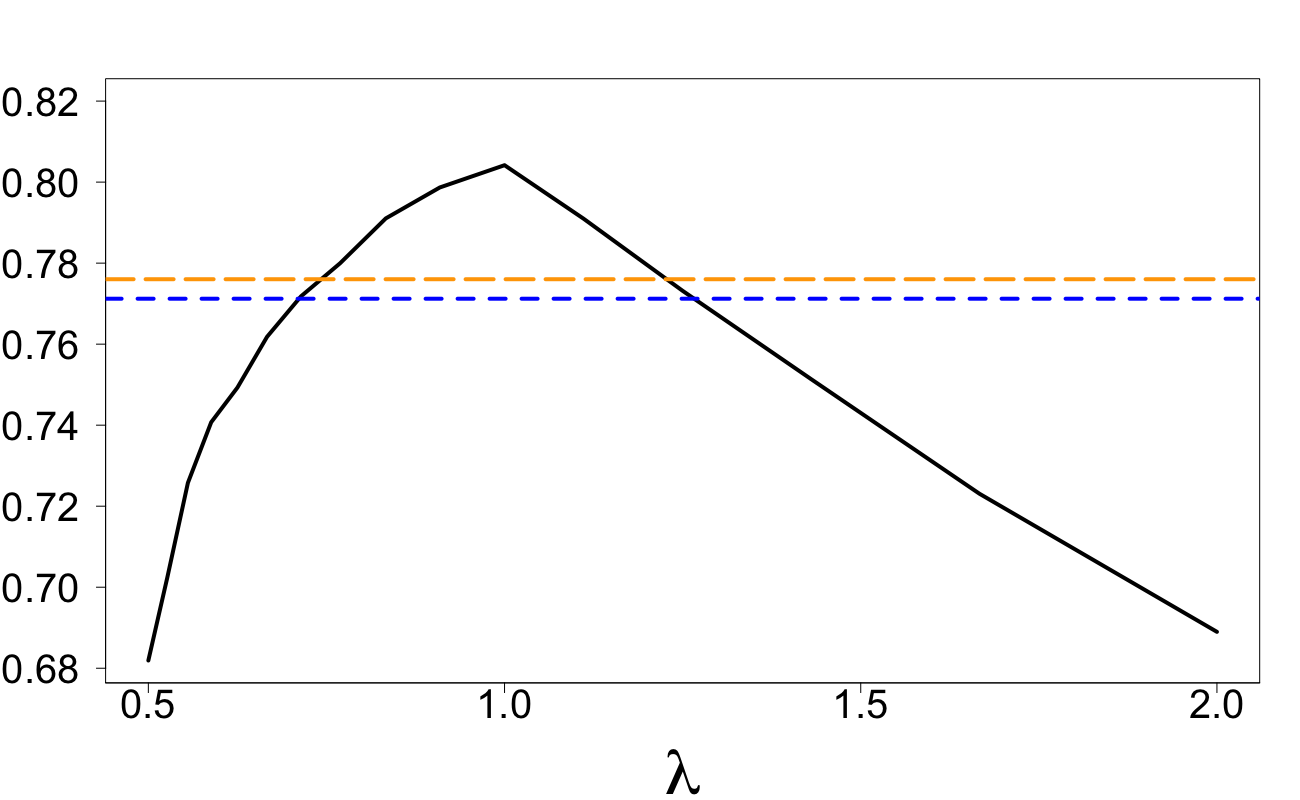}
 \end{center}
\caption{Power of three tests with $A=1$ and $l=10$ and ARL$=500$.  }
\label{power_over_lambda00}
\end{figure}

\begin{figure}
\begin{center}
\includegraphics[width=0.5\textwidth]{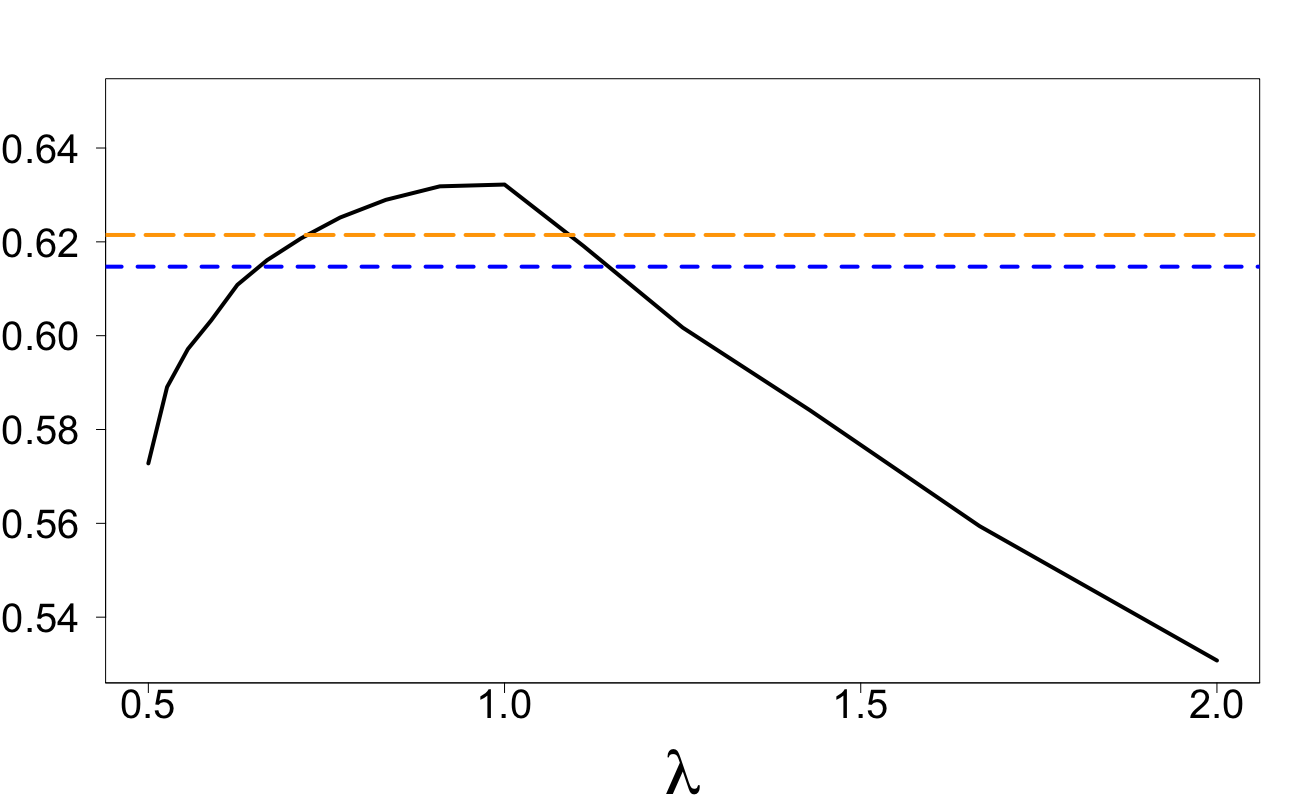}
 \end{center}
\caption{Power of three tests with $A=0.5$ and $l=20$ and ARL$=500$. }
\label{power_over_lambda11}
\end{figure}

In the first example shown in Figure~\ref{power_over_lambda00}, we have set $A=1$ and $l=10$. For the MOSUM test, we considered values of $L\in[5,20]$ to ensure $\lambda \in [0.5,2]$. For each $\lambda$, the values of ${\cal{P}}_S(H_1,A,L)$ can be accurately approximated using results in \cite{noonan2020power} or via Monte Carlo methods and are displayed with a solid black line. The dashed orange line depicts ${\cal{P}}_Z(H_2,A,5,20)$ which corresponds to prior knowledge that $l$ is between $[5,20]$. The shorter dashed blue line corresponds to ${\cal{P}}_{V}(H_{3},A)$ which has been obtained via Monte Carlo simulations. In  
Figure~\ref{power_over_lambda11}, we  set $A=0.5$ and $l=20$. For the MOSUM procedure, we consider values of $L\in[10,40]$ to ensure $\lambda \in [0.5,2]$. In this figure, the dashed orange line depicts ${\cal{P}}_Z(H_2,A,10,40)$ which corresponds to prior knowledge that $l$ is between $[10,40]$. The shorter dashed blue line corresponds to ${\cal{P}}_{V}(H_{3},A)$ obtained via Monte Carlo simulations.

From Figures~\ref{power_over_lambda00}-\ref{power_over_lambda11}, one can observe the advantage of knowing $l$ since the largest value of ${\cal{P}}_S(H_1,A,L)$ is the largest power of all three tests and is obtained for $\lambda=l/L=1$. In these figures, the values of $\lambda=l/L$ such that ${\cal{P}}_S(H_1,A,L)$ exceeds ${\cal{P}}_Z(H_2,A,5,20)$ (Figures~\ref{power_over_lambda00}) and ${\cal{P}}_Z(H_2,A,10,40)$ (Figures~\ref{power_over_lambda11}) shows the freedom in the choice of $L$ such that when $l$ is unknown, you still benefit over only assuming $l$ is bounded (similarly for CUSUM case when considering the dashed blue line). From these figures it is clear that unless you are very fortunate in choosing $L$ close to $l$ for the MOSUM test, you should use the generalised MOSUM test if $A$ is known. Unfortunately, there are no convenient analytic results for this test. Moreover, both the generalised MOSUM procedure and CUSUM procedures require the additional knowledge of $A$; this is not true for MOSUM. For the choice of parameters considered in both examples, the additional knowledge of a transient change leads to obvious benefits in power; those is seen by comparing the generalised MOSUM orange lines with the blue CUSUM lines.  Of course, ${\cal{P}}_Z(H_2,A,l_0,l_1)\rightarrow {\cal{P}}_S(H_1,A,l)$ as $l_0,l_1 \rightarrow l$.

\subsection{An application to real world data}\label{sec:real_world}

\begin{figure*}
\begin{center}
 \includegraphics[width=0.7\textwidth]{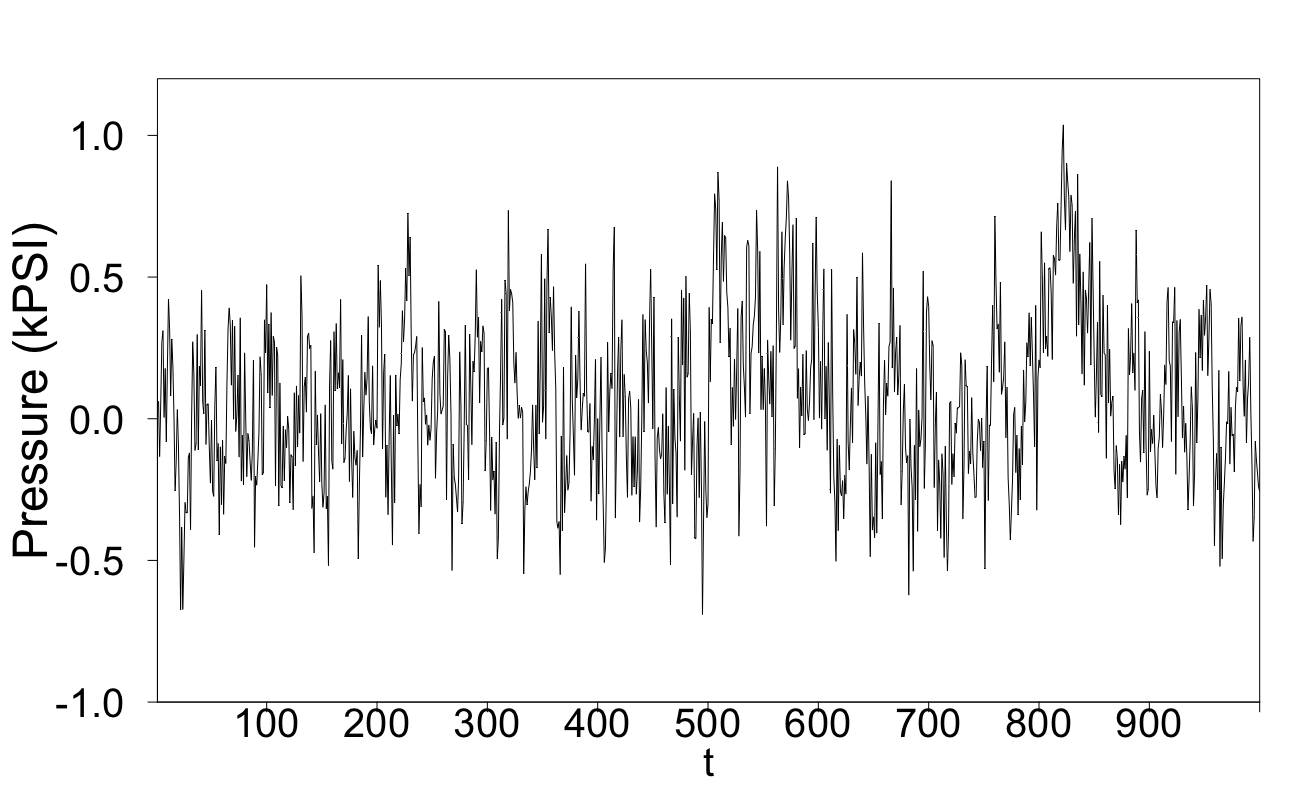}
\end{center}
\caption{Typical pressure data}
\label{example_data}
\end{figure*}

Hydrostatic pressure testing is important safety precaution for the Oil and Gas industry, see \cite{mcaleese2000operational}. Pressure testing is performed to confirm a pressure containing system is structurally sound and not leaking. Tests are performed by increasing the pressure in the system, expanding the pressure body, until the pressure reaches a pre-defined value typically equal to or larger than the body rated design pressure, then holding it there for a long enough time period to confirm there are no leaks, until eventually releasing the pressure. When performing tests offshore on floating Vessel/Drilling Rigs (Rig) this is complicated by the Rig’s movement due to the ocean waves, which introduce nearly sinusoidal fluctuations in pressure. Many of these tests are performed in real time and in parallel. Locating automatically when a test has been performed is essential for pressure analysis to determine if a leak is present and this is not obvious when noise is large. Typical example data is shown in Figure~\ref{example_data}.
When performing pressure tests, the hold periods can differ in length and amplitudes (pressure). 

\begin{figure}
\begin{center}
 \includegraphics[width=0.5\textwidth]{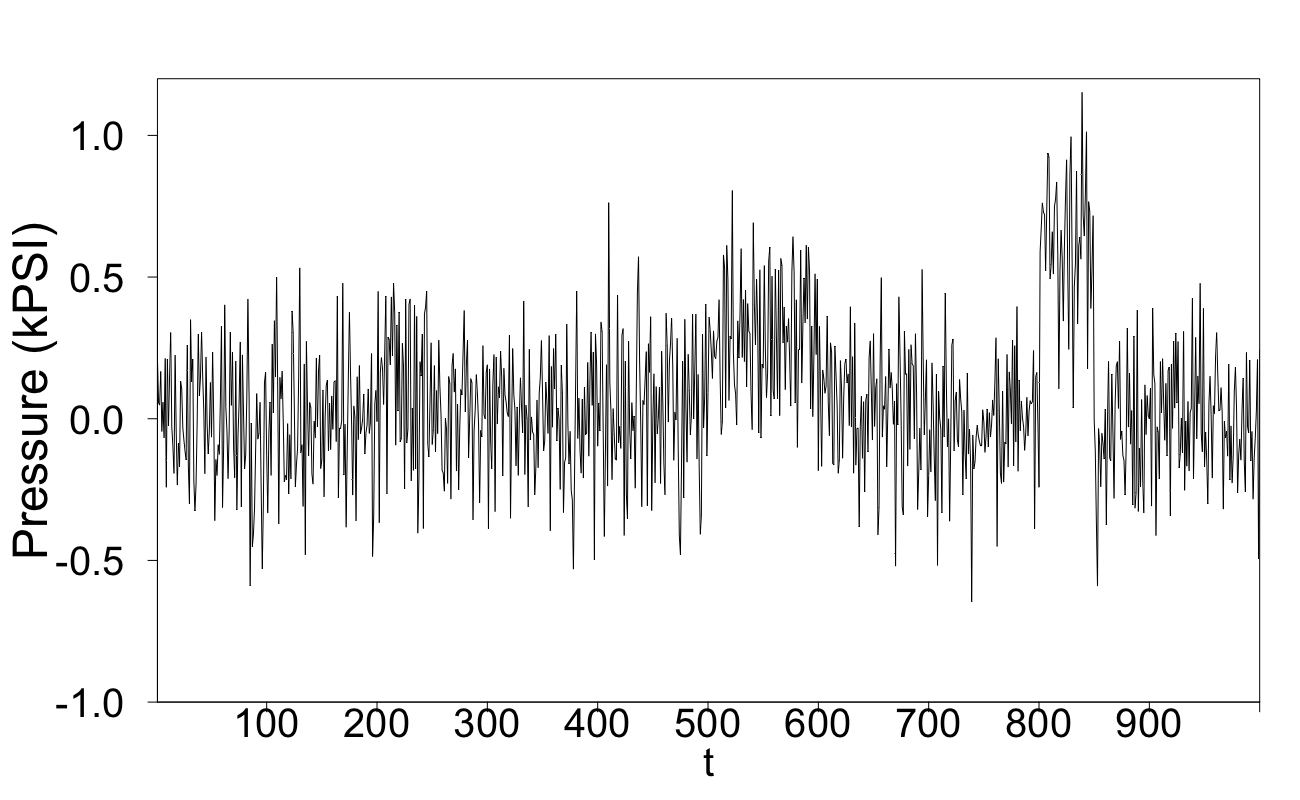}
 \end{center}
\caption{ Behaviour of $y_t$. }
\label{example_data2}
\end{figure}

\begin{figure}
\begin{center}
\includegraphics[width=0.5\textwidth]{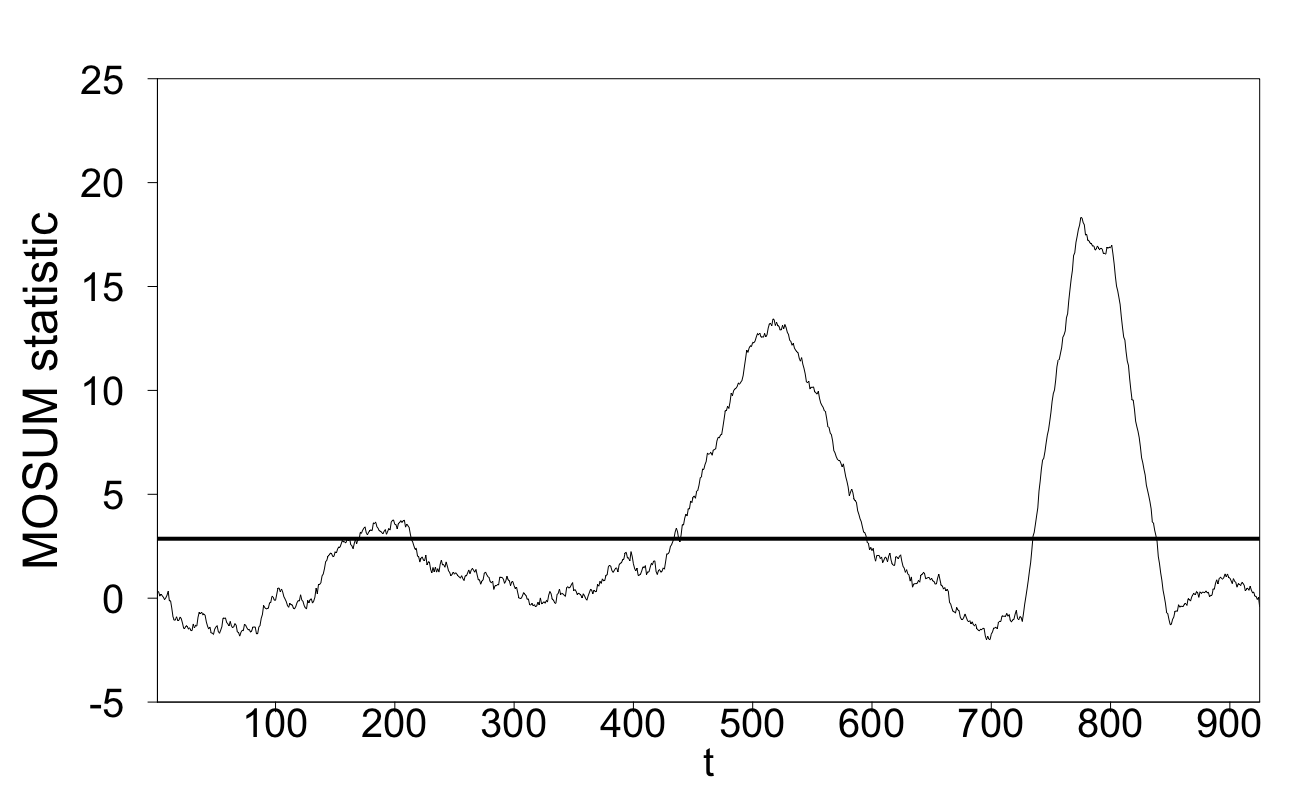}
 \end{center}
\caption{  MOSUM statistic based on Figure~\ref{example_data2} with $L=50$ and ARL$=5000$. }
\label{example_data3}
\end{figure}

\begin{figure}
\begin{center}
 \includegraphics[width=0.5\textwidth]{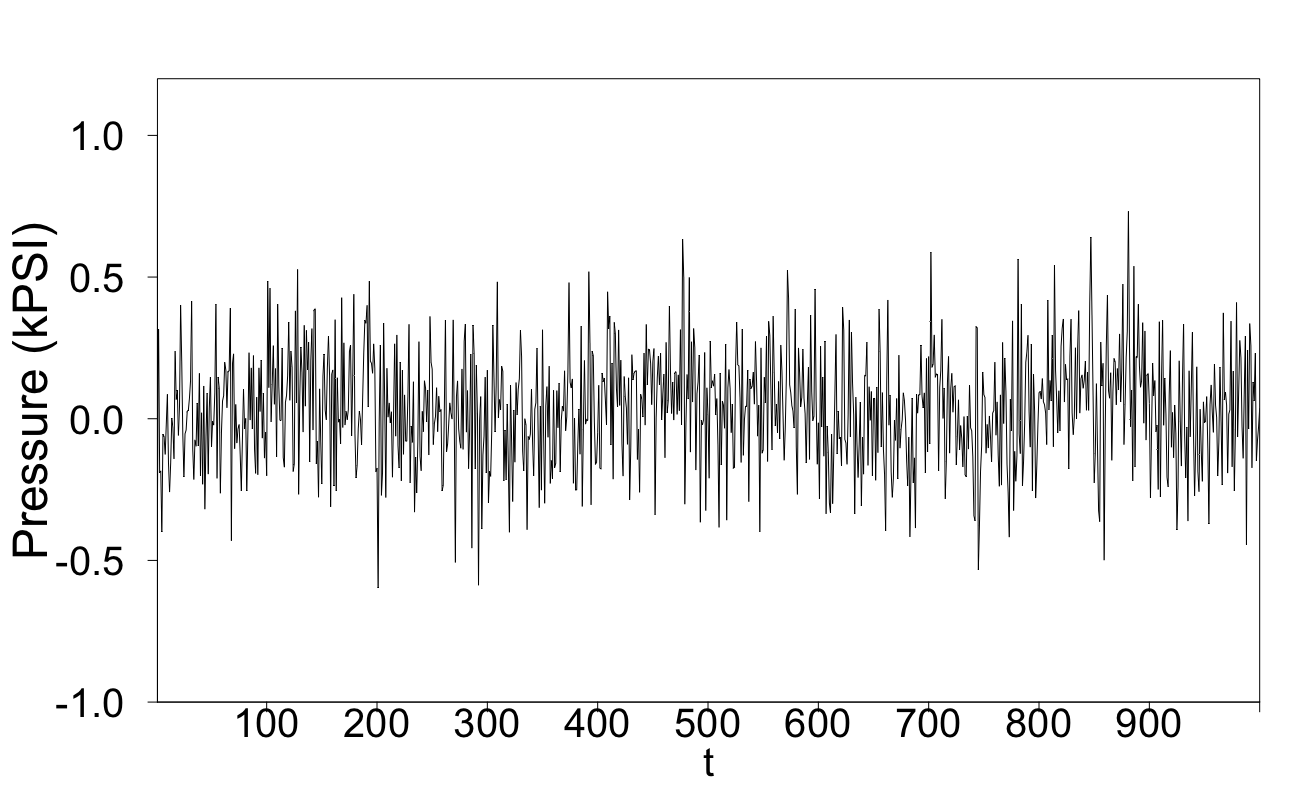}
 \end{center}
\caption{Behaviour of $y_t$. }
\label{example_data5}
\end{figure}

\begin{figure}[h]
\begin{center}
\includegraphics[width=0.5\textwidth]{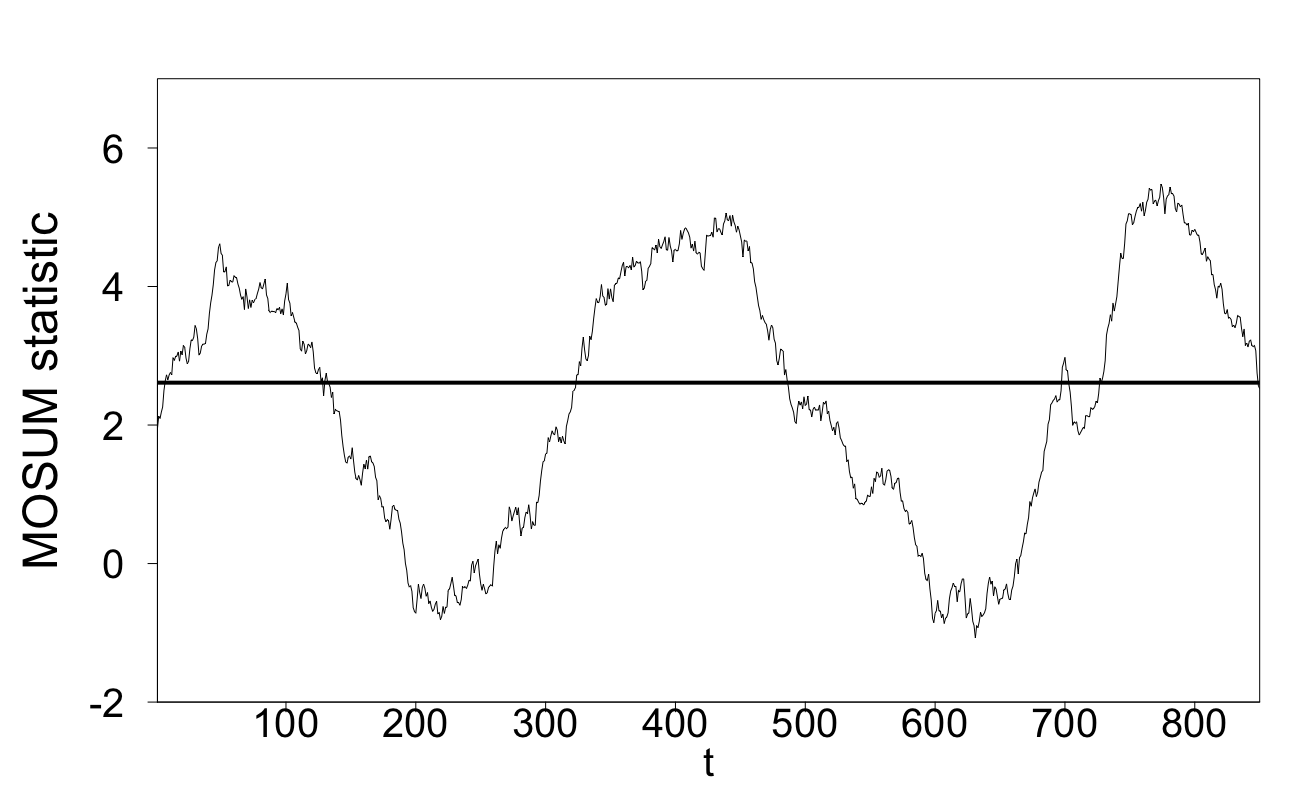}
 \end{center}
\caption{MOSUM statistic based on Figure~\ref{example_data5} with $L=150$ and ARL$=5000$. }
\label{example_data6}
\end{figure}

A sensible way of modelling the data under the null hypothesis of no pressure test could be
$
 z_t =s_t+y_t
$,
 where $s_t$ represents the signal introduced by the wave motion and $y_t$ can be modelled as i.i.d. $N(\mu,\sigma^2)$ and reflects the random noise that is present in the system. In most scenarios, there is significant pre-test data so $s_t, \mu$ and $\sigma$ can be estimated with great accuracy and therefore assumed known.  How to estimate $s_t$ or in general how to remove all main components of a signal leaving only noise can be performed using Singular Spectrum Analysis, see \cite{GNZ2001,golyandina2013singular}. When a pressure test begins, this can be reflected with a change in mean of the $y_t$; that is, under a pressure test $\mathbb{E}y_t =\mu+A $. The value of $A$ is often constant, but can differ between tests and is generally unknown. Each test can differ in duration but typical lengths vary between $l\in[50,100]$ units of time. One has to detect a transient change in mean of $y_t=z_t-s_t$. The behaviour of $z_t-s_t$ is shown in Figure~\ref{example_data2}.  In Figure~\ref{example_data3}, we depict the MOSUM statistic setting $L=75$. The horizontal line in this figure corresponds to the threshold required for an ARL of 5000. The MOSUM statistic indicates the location of three performed pressure tests and has the great advantage of not requiring knowledge $A$ when determining the ARL threshold unlike the generalised MOSUM and CUSUM procedures. A similar example is shown in Figures~\ref{example_data5}-\ref{example_data6}, where $L=150$ has been selected; three tests have been clearly located.

 \section*{Acknowledgements}
This work was partially supported by EPSRC grant EP/M024830.
 
 \bibliographystyle{unsrt}

\bibliography{changepoint}

%%  The body

%%  The bibliography

%%  If your bibliography is in BibTeX format, use the following setup:
%%  Style BST file for numbered citation:
%\bibliographystyle{imsart-number}
%%  or name-year citation
%\bibliographystyle{imsart-nameyear}
%%  Bibliography file (usually `*.bib')
%\bibliography{sii-sample}          
%%
%%  or include bibliography directly:
%\begin{thebibliography}{9}
%%  Use \bibitem{r1} or \bibitem[Surname(2010)]{r1} (for authoryear case)
%%  Put author names in \textsc{} command in order to use small caps font
%
%\bibitem{}
%\textsc{}
%
%\end{thebibliography}

\end{document}